\documentclass[a4paper]{article}
\usepackage{amsmath}
\usepackage{algorithm}
\usepackage{algpseudocode}
\usepackage{amsfonts}
\usepackage{authblk}
\usepackage[backend=bibtex]{biblatex}
\usepackage{booktabs}
\usepackage{hyperref}
\usepackage{cleveref}
\usepackage{bbold}
\usepackage{mathtools}
\usepackage{tikz}
\usepackage{caption}
\usepackage[acronym]{glossaries}
\usepackage{subcaption}
\usepackage{siunitx}
\usepackage{tabularx}
\usepackage{pgf,interval}
\usepackage[switch, modulo]{lineno}
\usepackage{tikz-dimline}
\usepackage{orcidlink}

\definecolor{waterColor}{RGB}{0,0,128}
\definecolor{groundColor}{RGB}{156,97,20}
\hypersetup{
	colorlinks,
	citecolor=waterColor,
	filecolor=black,
	linkcolor=black,
	urlcolor=waterColor
}
	
\usetikzlibrary{circuits.ee.IEC}
\usetikzlibrary{backgrounds,patterns,calc}
\usetikzlibrary{patterns.meta}

\newcommand{\dz}{\Delta z}
\newcommand{\dx}{\Delta x}
\newcommand{\dt}{\Delta t}
\newcommand{\expcvg}{\mathrm{CR}}
\newcommand{\orcid}[1]{\href{https://orcid.org/#1}{\textcolor[HTML]{A6CE39}{\aiOrcid}}}

\title{Convergence Properties of Iteratively Coupled Surface-Subsurface Models}
\author[*,1]{Valentina Schüller \orcidlink{0000-0002-4579-4217}}
\author[1]{Philipp Birken \orcidlink{0000-0002-6706-3634}}
\author[2]{Andreas Dedner \orcidlink{0000-0003-4387-9791}}
\affil[*]{Corresponding author, \href{mailto:valentina.schuller@math.lu.se}{valentina.schuller@math.lu.se}}
\affil[1]{Centre for Mathematical Sciences, Lund University, Lund, Sweden}
\affil[2]{Department of Mathematics, University of Warwick, Coventry, United Kingdom}
\date{}

\begin{document}

\maketitle

\begin{abstract}
  Surface-subsurface flow models for hydrological applications solve a coupled multiphysics problem.
  This usually consists of some form of the Richards and shallow water equations.
  A typical setup couples these two nonlinear partial differential equations in a partitioned approach via boundary conditions.
  Full interaction between the subsolvers is ensured by an iterative coupling procedure.
  This can be accelerated using relaxation.

  In this paper, we apply continuous and fully discrete linear analysis techniques to study an idealized, linear, 1D-0D version of a surface-subsurface model.
  These result in explicit expressions for the convergence factor and an optimal relaxation parameter, depending on material and discretization parameters.
  We test our analysis results numerically for fully nonlinear 2D-1D experiments based on existing benchmark problems.
  The linear analysis can explain fast convergence of iterations observed in practice for different materials and test cases, even though we are not able to capture various nonlinear effects.
\end{abstract}

\paragraph{Competing Interests} The authors declare no competing interests.


\section{Introduction}\label{sec:introduction}

We consider a coupled surface-subsurface flow model. 
This is motivated by the desire to better understand the terrestrial water cycle and to simulate effects of water management and flooding on the water table and aquifer \parencite{keyes.etal_2013,bastian.etal_2012}. 
A common choice is to use a $d$-dimensional Richards equation for the subsurface flow and some variant of the $d-1$-dimensional shallow water equations for the surface flow. 
We thus have two coupled nonlinear partial differential equations (PDEs). 

We follow a so-called partitioned approach. 
There one uses separate solvers, with potentially different discretizations, for the two PDEs. 
The models then interact by prescribing boundary conditions for each other. Examples of coupled surface-subsurface models that make use of such a strategy are given in \textcite{furman_2008,maxwell.etal_2014}.
Specifically, the shallow water equation supplies the water height as a Dirichlet condition for the capillary head, and the Richards equation supplies a flux in the form of the vertical velocity. 
The latter acts as a source term in the lower dimensional surface flux equation.
As in domain decomposition methods, this is performed in an iterative manner. 

While this strategy is used successfully in practice, an analysis of this scheme is to our knowledge missing for coupled surface-subsurface flows. 
Additionally, convergence can be accelerated using a relaxation step, but this is not used so far. 
In this article, our goal is to provide a convergence analysis of the coupling iteration and determine optimal relaxation parameters.

Such an analysis has been performed previously for linear coupled problems, in particular for heterogeneous heat equations.
This has been done on the fully continuous level \parencite{gander.etal_2016}, discrete in time \parencite{henshaw.chand_2009}, discrete in space \parencite{janssen.vandewalle_1996a}, or fully discrete \parencite{monge.birken_2018}. 
The convergence results depend on the material parameters and, in the latter three cases, on the discretization.
The fully discrete results show that the linear analysis is able to predict the quantitative and qualitative behavior of the convergence factor for the nonlinear case.

We apply this type of linear analysis to the nonlinear problem of coupled surface-subsurface flow, where in particular for Richards' equation material parameters vary in size by many orders of magnitude.
To this end, we consider a linearized 1D-0D model.
We then focus on the fully discrete approach from \textcite{monge.birken_2018}, which is able to show the effect of the different mesh widths on convergence. 
To get even more insight and to prove well-posedness, we also apply the continuous analysis from \textcite{gander.etal_2016}. 
Both result in explicit expressions for the convergence factor of the iterations as well as optimal relaxation parameters.

We then look at the fully nonlinear coupled 2D-1D problem and carry out numerical experiments with realistic material parameters based on existing benchmark problems.
These allow us to measure speed of convergence in practice and compare with the results predicted by our linear analysis.
The two subproblems are solved with the PDE software framework DUNE \parencite{dedner.etal_2010} and coupled using the coupling library preCICE \parencite{chourdakis.etal_2022}.

In \Cref{sec:governing_equations}, we first present the nonlinear coupled problem and then in \Cref{sec:discretization}, the discretizations and the coupling method. 
The analysis of the method is presented in \Cref{sec:linear_analysis} and in \Cref{sec:numerical_results}, we show numerical experiments to validate the linear analysis and compare with the fully nonlinear code. 
Conclusions are presented in \Cref{sec:conclusion}.

\section{Governing Equations}\label{sec:governing_equations}

We study two-dimensional coupled subsurface-surface flow, i.e., the subsurface model is defined on a domain \(\Omega\subset\mathbb{R}^2\).
The surface model lives on the one-dimensional upper boundary of \(\Omega\).
Our model is similar to the ones presented in \textcite{bastian.etal_2012,kollet.maxwell_2006}, albeit in 2D instead of 3D.

\subsection{Subsurface Flow}\label{subsec:subsurface_model}

\begin{figure}
    \centering
    \begin{tikzpicture}
    \tikzstyle{doublefullarrow}=[<->, >=triangle 60, very thin]
    \tikzset{>=latex}
    \draw[draw=black, pattern=dots, pattern color=groundColor] (0.0, 0.0) rectangle ++(3.0,2.0);
    \draw[waterColor, very thick] (0.0, 2.03) -- (3.0, 2.03);
    \node[waterColor, anchor=south] (Gamma) at (1.5, 2.0) {\footnotesize $\Gamma$};
    \node[black] at (1.5, 1) {\large $\Omega$};
    \node[anchor=south] (dOmega) at (3.5, 1.5) {\footnotesize $\partial\Omega$};
    \draw[black, -] (dOmega) -- (3, 1);
    \node[shape=coordinate] (CSorigin) at (0.2, 0.2) {};
    \node[] (CSz) at (0.2, 0.8) {\footnotesize $z$};
    \node[] (CSx) at (0.8, 0.2) {\footnotesize $x$};
    \draw[black, ->] (CSorigin) -- (CSz);
    \draw[black, ->] (CSorigin) -- (CSx);
\end{tikzpicture}
    \caption{The spatial domain in the 2D-1D coupling problem.\label{fig:domain_2d}}
\end{figure}
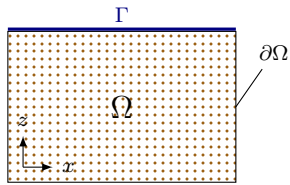

Subsurface flow is usually modelled with some form of the Richards equation \parencite{furman_2008,maxwell.etal_2014,keyes.etal_2013}.
Richards' equation can be expressed in terms of the volumetric soil water content \(\theta(t,\boldsymbol{x})\) or the capillary head \(\psi(t,\boldsymbol{x})\), where \(t\) denotes time and \(\boldsymbol{x}=(x,z)^T\) the position.
As in \textcite{bastian.etal_2012,list.radu_2016}, we look at the Richards equation in so-called mixed form,
\begin{equation}\label{eq:richards_2d}
    \begin{aligned}
        \partial_t \theta(\psi) 
        + \nabla\cdot\underbrace{\left(-K(\psi)\nabla(\psi+z) \right)}_{\eqqcolon \boldsymbol{v}(\psi)} 
        &= 0 \quad \text{on } (0,T] \times \Omega, \\
        \boldsymbol{v}\left(\psi(t,\boldsymbol{x})\right) \cdot \boldsymbol{n} &= 0 \quad \text{on } [0,T] \times \partial\Omega \setminus \Gamma,\\
        \psi(t,\boldsymbol{x}) &= \psi_\Gamma(t,x) \quad \text{on } [0,T] \times \Gamma,\\
        \psi(0, \boldsymbol{x}) &= \psi_0(\boldsymbol{x}) \quad \text{on } \Omega.
    \end{aligned}
\end{equation}
Here, \(\theta(\psi)\) is the soil water content, \(K(\psi) = K(\theta(\psi))\) is the unsaturated hydraulic conductivity, \(\boldsymbol{v}(\psi)\) is the volumetric flux, and \(\boldsymbol{n}\) denotes the outward-facing unit normal vector on \(\partial\Omega\).
We consider the spatial domain \(\Omega=(0,L_x)\times(0,L_z)\) with boundary \(\partial\Omega\), cf. \Cref{fig:domain_2d}.
The subset \(\Gamma=\{\boldsymbol{x} | z = L_z)\}\subset \partial\Omega\) is the upper boundary of the subsurface domain.

Both \(\theta\) and \(K\) nonlinearly depend on the capillary head \(\psi\) via empirical constitutive relations, for which we use the van Genuchten--Mualem model \parencite{vangenuchten_1980}:
\begin{equation}\label{eq:genuchten}
    \begin{aligned}
        \theta(\psi) &= \begin{cases}
            \theta_R + (\theta_S - \theta_R) \left(\frac{1}{1+(\alpha |\psi|)^n_G}\right)^{\frac{n_G-1}{n_G}}, & \psi \leq 0, \\
            \theta_S, & \psi > 0,
        \end{cases}\\
        K(\psi) &= \begin{cases}
            K_S \sqrt{\theta(\psi)} \left[1-\left(1-\theta(\psi)^{\frac{n_G}{n_G-1}}\right)^\frac{n_G-1}{n_G}\right]^2, & \psi \leq 0, \\
            K_S, & \psi > 0.
        \end{cases}
    \end{aligned}
\end{equation}
Therein, \(\theta_S\) and \(\theta_R\) are the saturated and residual water content, \(K_S\) is the saturated hydraulic conductivity, and \(\alpha\) and \(n_G\) are soil parameters related to the air entry suction and pore-size distribution, respectively.
We will later also make use of the hydraulic capacity \(c(\psi)=\partial\theta/\partial\psi\), given by
\begin{equation}\label{eq:capacity_definition}
    c(\psi) = \begin{cases}
        \alpha (\theta_S - \theta_R)(n_G - 1)(\alpha |\psi|)^{n_G - 1}
        (1 + (\alpha |\psi|)^{n_G}) ^{(1 / n_G - 2)}, & \psi \leq 0, \\
            0, & \psi > 0.
    \end{cases}
\end{equation}

For simplicity, we assume no-flux boundary conditions at the left, right, and lower boundary here.
The extension to general boundary conditions is straightforward \parencite{list.radu_2016}.

Using the method of weighted residuals and integration by parts, we obtain the following weak form of \Cref{eq:richards_2d}:
Find \(\psi\in H^1(\Omega)\) such that
\begin{equation*}
    \left\langle\partial_t \theta(\psi), \phi \right\rangle 
    - \left\langle \boldsymbol{v}(\psi), \nabla \phi \right\rangle 
    = f_\Gamma
\end{equation*}
holds for all \(\phi\in H^1(\Omega)\).
The right-hand side \(f_\Gamma\) contains the contribution from the Dirichlet boundary \(\Gamma\):
\begin{equation}\label{eq:richards_fgamma}
    f_\Gamma = -\left\langle\partial_t \theta\left(\psi_\Gamma\right), \phi \right\rangle 
    + \left\langle \boldsymbol{v}\left(\psi_\Gamma\right), \nabla \phi \right\rangle. 
\end{equation}

\subsection{Surface Flow}\label{subsec:surface_model}

We consider the special case where the uppermost layer of the ground is fully saturated and surface water is present. If this is not the case, the water height is zero and the surface flow model is inactive \parencite{sulis.etal_2010}. 

In a first approach, we model surface flow with the one-dimensional shallow water equations (SWE) for the water height \(h(t,x)\) and velocity \(u(t,x)\), cf. \textcite{bastian.etal_2012}:
\begin{equation}\label{eq:swe}
    \begin{aligned}
        \partial_t \begin{pmatrix}h \\ hu\end{pmatrix} 
        + \partial_x \begin{pmatrix} hu \\ hu^2 + \frac{1}{2} gh^2 \end{pmatrix}
        &= \begin{pmatrix} s(t,x) + r(t,x) \\ 0 \end{pmatrix} \quad \text{on } (0,T] \times \Gamma,\\ 
        \begin{pmatrix}h \\ hu\end{pmatrix}(0, x) &= \begin{pmatrix}h_0 \\ h_0u_0\end{pmatrix} (x) \quad \text{on } \Gamma.
    \end{aligned}
\end{equation}
Here, \(g=9.81 \unit{m.s^{-2}}\) denotes acceleration by gravity, \(s(t,x)\) is a source and sink term for the water height, and $r(t,x)$ denotes the rainfall rate.
In our model problem, we do not consider additional terms due to bathymetry or friction. 
This model is suitable for small-scale surface flow.

In a second approach, we follow a similar implementation as large-scale coupled hydrology models such as ParFlow \parencite{kollet.maxwell_2006} or CATHY \parencite{camporese.etal_2010} and make use of the kinematic wave approximation.
The model reduces to a single equation for conservation of mass:
\begin{equation}\label{eq:maxwell}
	\begin{aligned}
        \partial_t h + \partial_x (hu) &= s(t,x) + r(t, x) \quad \text{on } (0,T] \times \Gamma,\\ 
        h(0,x) &= h_0(x) \quad \text{on } \Gamma,
    \end{aligned}
\end{equation}
where the velocity \(u\) is determined from Manning's equation,
\begin{equation*}
	u = \frac{\sqrt{S_{f}}}{n_M} h^{2/3}.
\end{equation*}
Therein, $S_f$ is the bed slope and $n_M$ is Manning's coefficient.

We can write both equations as a one-dimensional balance law with the shorthands \(\boldsymbol{q}={(h,hu)}^T\), \(\boldsymbol{b}={(s + r, 0)}^T\) and $\boldsymbol{q} = h$, $\boldsymbol{b} = s + r$ for \Cref{eq:swe,eq:maxwell}, respectively.
Assuming homogeneous Neumann conditions at the boundary to be consistent with the subsurface model, we arrive at:
\begin{equation}\label{eq:swe_balance_law}
    \begin{aligned}
        \partial_t \boldsymbol{q}(t,x) + \partial_x f\left(\boldsymbol{q}(t,x)\right) &= \boldsymbol{b}(t,x) \quad \text{on } (0,T] \times \Gamma,\\
        \partial_x f\left(\boldsymbol{q}(t,0)\right) = \partial_x f\left(\boldsymbol{q}(t,L_x)\right) &= 0 \quad \text{on } [0,T],\\
        \boldsymbol{q}(0, x) &= \boldsymbol{q}_0(x) \quad \text{on } \Gamma.
    \end{aligned}
\end{equation}

\subsection{Interface Conditions and Coupling Iteration}\label{subsec:cont_coupling}

Surface-subsurface models can be coupled at the interface \(\Gamma\) either using boundary conditions that enforce continuity of pressure and mass conservation \parencite{bastian.etal_2012,kollet.maxwell_2006}, or alternatively using a first-order exchange flux \parencite{vanderkwaak.loague_2001,panday.huyakorn_2004}.
We follow the first approach, giving for \(t\in(0,T]\):
\begin{subequations}\label{eq:nonlinear_ibc}
    \begin{align}
        \psi_\Gamma(t,x) &= h(t, x), \label{eq:ibc_dirichlet} \\
        s(t, x) &= \boldsymbol{v}(\psi(t, x, L_z)) \cdot \boldsymbol{n} \eqqcolon v_s(t,x). \label{eq:ibc_neumann}
    \end{align}
\end{subequations}
We focus on the situation with nonzero water height. 
Coupled hydrology models need to account for cases with \(h=0\) as well, this is not the focus of this paper.
In that case, the interface conditions \eqref{eq:nonlinear_ibc} no longer hold.\footnote{\(\psi(t,0)<0\), but \(h=0\), negative water heights work with neither of the two surface flow models.}
Examples of how this is tackled in existing codes are given in \textcite{sulis.etal_2010,maxwell.etal_2014}.
There exists no known analytical solution to this coupled problem \parencite{maxwell.etal_2014,bastian.etal_2012}.

We solve the coupled problem using partitioned, sequential coupling iterations with relaxation.
Given an initial guess for the water height \(h^0(t,x)\), compute sequentially in each iteration \(k=1,2,\dots\) the approximations for \(\psi^k\) followed by \(\boldsymbol{\tilde{q}}^k\):
\begin{subequations}\label{eq:cont_cpl_iterations}
\begin{equation}
    \begin{aligned}
        \partial_t \theta(\psi^k) 
        + \nabla\cdot\boldsymbol{v}(\psi^k) &= 0 \quad \text{on } (0,T] \times \Omega, \\
        \boldsymbol{v}\left(\psi^k(t,\boldsymbol{x})\right) \cdot \boldsymbol{n} &= 0 \quad \text{on } [0,T]\times\partial\Omega \setminus\Gamma, \\
        \psi^k(0,\boldsymbol{x}) &= \psi_0(\boldsymbol{x}) \quad\text{on }\Omega, \\
        \psi^k(t,\boldsymbol{x}) &= h^{k-1}(t,x) \quad \text{on } [0,T]\times \Gamma,
    \end{aligned}
\end{equation}
\begin{equation}
    \begin{aligned}
        \partial_t \boldsymbol{\tilde{q}}^k + \partial_x f(\boldsymbol{\tilde{q}}^k) &= \boldsymbol{b}^k \quad \text{on } (0,T]\times \Gamma, \\
        \partial_x f\left(\boldsymbol{\tilde{q}}(t,0)\right) = \partial_x f\left(\boldsymbol{\tilde{q}}(t,L_x)\right) &= 0 \quad \text{on } [0,T],\\
        \boldsymbol{\tilde{q}}(0, x) &= \boldsymbol{q}_0(x) \quad \text{on } \Gamma, \\
        s^k(t,x) &= \boldsymbol{v}(\psi^k(t,x,L_z)) \cdot \boldsymbol{n} \quad \text{on } [0,T]\times \Gamma.
    \end{aligned}
\end{equation}
Then compute the new water height \(h^k(t,x)\) using
\begin{equation}
    h^k(t,x) = \omega \tilde{h}^k + (1-\omega) h^{k-1},
\end{equation}
\end{subequations}
with the relaxation parameter \(\omega \in (0,1]\).
The main goal of this paper is to study the convergence properties of this coupling iteration and how they relate to the choice of \(\omega\).

This type of equation coupling is reminiscent of Dirichlet-Neumann waveform relaxation (DNWR) in domain decomposition or fluid-structure interaction:
The surface solver provides a Dirichlet condition for the subsurface model, while the subsurface solver sends back a (nonlinear) flux.
However, the flux enters the surface flow model as a source term, not a boundary condition, which is a fundamental difference to the classical DNWR algorithm.

\section{Discretization of the Coupling Iterations}\label{sec:discretization}

This section presents a fully discrete version of the coupling iteration \eqref{eq:cont_cpl_iterations}.
We follow the method of lines, discretizing first in space, then in time.

\subsection{Discretization of the Subsurface Model}

For the subsurface model, we use linear continuous finite elements in space and the implicit Euler method in time.\footnote{Using an implicit first-order time discretization for Richards' equation is a common choice \parencite{farthing.ogden_2017}.}
We solve Richards' equation on a Cartesian grid with \(M_x \times M_z\) elements of size \(\dx = L_x/M_x\) and \(\dz = L_z/M_z\).
Furthermore, we introduce uniformly spaced time steps \(t^n\), with \( \Delta t = T/N \).
Then we find approximate solutions \( \psi^{n,k}(\boldsymbol{x}) \in \mathcal{V}^h \) to \(\psi^k(t^n, \boldsymbol{x})\), where \(\mathcal{V}^h\) is a finite-dimensional subspace of \(H^1(\Omega)\). 

We use a linear nodal basis for \(\psi^{n,k}\) and \(\phi\) in our implementation and arrive at the following nonlinear system:
Given \(\psi^{n-1} \in \mathcal{V}^h\), find \(\psi^{n,k} \in \mathcal{V}^h\) such that
\begin{equation}\label{eq:discr_richards_wf}
    \left\langle \theta\left(\psi^{n,k}\right) - \theta\left(\psi^{n-1}\right), \phi \right\rangle 
    - \dt \left\langle \boldsymbol{v}(\psi^{n,k})\cdot\boldsymbol{n}, \nabla \phi \right\rangle 
    = f_\Gamma^{n,k} \quad \forall \phi \in \mathcal{V}^h.
\end{equation}
The initial condition \(\psi^0\) is obtained by suitably approximating \(\psi_0(\boldsymbol{x})\) in \(\mathcal{V}^h\).
The nonlinear system for $\psi^{n,k}(\boldsymbol{x})\approx \psi^k(t^n, \boldsymbol{x})$ is solved using Newton's method.

\subsection{Discretization of the Surface Model}

We discretize the surface model using a cell-centered finite volume scheme in space and the implicit Euler method in time.\footnote{Note that using the implicit Euler method for the surface flow solver is not uncommon in coupled hydrological models, cf. \cite{kollet.maxwell_2006}.}
Without loss of generality, we assume an equidistant mesh for the surface flow domain \(\Gamma\).
This consists of \(M_x\) cells \(\mathcal{C}_l = [x_{l-1}, x_l]\), with \(x_l = l\dx\) and \(l=1,\dots,M_x\).
We work with the balance law formulation given in \Cref{eq:swe_balance_law}.
Defining
\begin{equation*}
    \boldsymbol{\tilde{q}}_l^k(t) \coloneqq \frac{1}{\dx} \int_{\mathcal{C}_l} \boldsymbol{\tilde{q}}^k(t,x) \, dx,
\end{equation*}
and introducing a numerical flux \(F_{l}^k(t)\), we obtain
\begin{equation}
    \frac{d}{dt} \boldsymbol{\tilde{q}}^k_l(t) = \frac{1}{\dx} \int_{\mathcal{C}_l} \boldsymbol{b}^k(t,x) \, dx - \frac{1}{\dx} \left(F_{l}^k(t) - F_{l-1}^k(t)\right).
\end{equation}

In our tests, we use the local Lax-Friedrichs method to define the numerical flux.
Integration in time and applying the implicit Euler method yields
\begin{equation}\label{eq:discr_swe}
    \boldsymbol{\tilde{q}}_l^{n,k} = \boldsymbol{q}_l^{n-1} - \frac{\dt}{\dx} \left(F_{l}^{n,k} - F_{l-1}^{n,k}\right) + \frac{\dt}{\dx} \int_{\mathcal{C}_l} \boldsymbol{b}^{n,k}(x) \, dx.
\end{equation}
Note that we use the same time step size \(\dt\) and horizontal grid size \(\dx\) for surface and subsurface flow.

\subsection{Fully Discrete Partitioned Coupling Iteration}\label{subsec:discr_cpl_iteration}

To finalize the discretization of the coupled problem, we need to address the interface boundary conditions \eqref{eq:nonlinear_ibc}.
Both solvers use the same time integration method and time step size. However, we need to provide a map from the cell-centered FV grid (surface grid) to the nodal values of the subsurface grid and vice versa.

The surface grid consists of $M_x$ cells, centered at $\frac{\Delta x}{2},\frac{3\Delta x}{2},\dots,L_x-\frac{\Delta x}{2}$, whereas the subsurface grid consists of $M_x+1$ nodes at $0,\dx,\dots,L_x$.
To map from the surface to the subsurface solver, we apply the rule that a nodal value is the average of the FV cells surrounding it:
\begin{equation}\label{eq:ibc_gnd}
    \psi_\Gamma^{n,k}(x_l) = 
    \begin{cases}
        \frac{1}{2} \left(h_{l}^{n,k-1} + h_{l+1}^{n,k-1}\right) & l=1,\dots,M_x-1 \\
        h_1^{n,k-1} & l=0 \\
        h_{M_x}^{n,k-1} & l=M_x.
    \end{cases}
\end{equation}

Vice versa, the subsurface solver provides the integral of the surface flux:
\begin{equation}\label{eq:ibc_sfc}
    \int_{\mathcal{C}_l} s^{n,k}(x) \, dx = \int_{x_{l-1}}^{x_l} \boldsymbol{v}\left(\psi^{n,k}(x,L_z)\right) \cdot \boldsymbol{n}\, dx = v_s^{n,k}\left(x_{l-1} + \frac{\dx}{2}\right) \dx,
\end{equation}
for each of the \(M_x\) cells \(\mathcal{C}_l\).
The second equality holds since we use linear finite elements.

The partitioned, sequential coupling iteration procedure in the fully discrete case takes the following form:
Given an initial guess for the water height \(h_l^{n,0}\), first solve \Cref{eq:discr_richards_wf} with the boundary condition \eqref{eq:ibc_gnd}.
Then, update the surface solver's source term by \Cref{eq:ibc_sfc} and solve \Cref{eq:discr_swe} to obtain \(\boldsymbol{\tilde{q}}_l^{n,k}\).
Finally, apply a relaxation step to update the Dirichlet boundary condition for the next iteration: 
\begin{equation}
    h_l^{n,k} = \omega \tilde{h}_l^{n,k} + (1-\omega) h_l^{n,k-1},
\end{equation}
with the relaxation parameter \(\omega \in (0,1]\).
We terminate the iteration once the residual for the water height falls below a given tolerance:
\begin{equation}\label{eq:conv_crit}
    \mathrm{res}^{n,k} \coloneqq \left\| \tilde{h}^{n,k} - h^{n,k-1} \right\|_2<\mathrm{TOL}. 
\end{equation}

\section{Linear Analysis}\label{sec:linear_analysis}

In this section, we analyze a linearized, one-dimensional version of the coupling problem described in \Cref{sec:governing_equations}.
We use techniques from \textcite{gander.etal_2016} and \textcite{monge.birken_2018} to obtain a convergence rate and an optimal relaxation parameter for the continuous and fully discrete coupling iteration, respectively.

\subsection{One-Dimensional Idealized Model}

We consider a single column of the coupled system, assuming that the convergence behavior of the iterations is driven by vertical processes.
Richards' equation reduces to a one-dimensional PDE whereas both surface flow models reduce to the same ordinary differential equation (ODE) for the water height:
\begin{equation*}
	\begin{gathered}
		\partial_t \theta(\psi) + \partial_z \left(K(\psi)\partial_z\left(\psi(t,z)+ z\right)\right) = r(t), \\
		\frac{dh}{dt} = v_s(t).
	\end{gathered}
\end{equation*}

Momentum conservation in \Cref{eq:swe} trivially becomes \(\dot{hu} = 0\).
Since it is no longer coupled to the continuity equation, we can disregard it in the convergence analysis.
To linearize Richards' equation, we write it in fully head-based form via
\begin{equation*}
	\frac{\partial \theta \left(\psi(t,z)\right)}{\partial t} = \underbrace{\frac{\partial\theta(\psi)}{\partial\psi}}_{\eqqcolon c(\psi)} \frac{\partial\psi(t,z)}{\partial t},
\end{equation*}
and assume that \(c,K>0\) are positive constants independent of \(\psi\).\footnote{Positivity is fulfilled in practice for \(K\), cf. \Cref{eq:genuchten}.
However, \(c(\psi)=0\) in saturated regions where \(\psi>0\), i.e., the equation becomes elliptic.
This degenerate nature of Richards' equation makes it hard to study analytically and numerically \parencite{farthing.ogden_2017}.}
We furthermore drop the rainfall rate, \(r(t)=0\), since it does not affect the convergence rate in the linear setting.
As in \Cref{subsec:cont_coupling}, we write down the partitioned coupling iteration including relaxation for this continuous model problem.

Given an initial guess \(h^0(t)\), compute for each \(k=1,2,\dots\) the approximations \(\psi^k\) and \(\tilde{h}^k\) sequentially according to:
\begin{subequations}\label{eq:linear_continuous_cpl}
\begin{equation}\label{eq:linear_cont_richards}
	\begin{aligned}
		c \frac{\partial\psi^k}{\partial t} - K\frac{\partial^2}{\partial z^2} \psi^k & = 0 \quad \text{on }(0, T]\times(-L, 0), \\
		\psi^k(0,z) & = \psi_0(z) \quad \text{on }[-L, 0],  \\
		\psi^k(t,-L) &= 0 \quad \text{on } [0, T], \\
		\psi^k(t,0) &= h^{k-1}(t) \quad \text{on } [0, T],
	\end{aligned}
\end{equation}
\begin{equation}
	\begin{aligned}\label{eq:linear_cont_swe}
		\frac{d\tilde{h}^k}{dt} & = v_s^k(t)\quad \text{on } (0, T], \\
		\tilde{h}^k(0) & = h_0.
	\end{aligned}
\end{equation}
The boundary data for the subsurface solver in the next iteration is updated by the relaxation step:
\begin{equation}
	h^k(t) = \omega \tilde{h}^k(t) + (1-\omega) h^{k-1}(t).
\end{equation}
\end{subequations}

Note that we use homogeneous Dirichlet boundary conditions for Richards' equation here, instead of no-flux boundary conditions as in the nonlinear problem.
This simplifies the fully discrete analysis and is justified if we assume that the bottom boundary is far away from the interface.\footnote{We have verified that the lower boundary condition does not meaningfully affect the convergence rate in our code. Switching from a Dirichlet-type to a Neumann-type boundary condition leads to changes in the observed convergence rate below \qty{10}{\percent} in magnitude.}

\Cref{eq:linear_continuous_cpl} corresponds to a homogeneous heat equation, sequentially coupled to a linear ODE.
The underlying coupled problem is well-posed under the assumption that $c,K,L>0$, see \Cref{subsec:well_posedness}.

\subsection{Continuous Analysis}\label{subsec:cont_cpl_analysis}

One can determine the convergence factor of the continuous iterations, as well as an optimal relaxation parameter, with Fourier-based techniques presented in \textcite{gander.etal_2016}.
We use the Laplace transform for time-dependent functions
\begin{equation}\label{eq:laplace_transform}
    \hat{f}(s) \coloneqq \mathcal{L}\{f\}(s) = \int_0^\infty f(t) e^{-st} \,dt,
\end{equation}
where \(s \in \mathbb{C}\) is the Laplace variable.
We can then state the sequential coupling iteration including relaxation for the Laplace-transformed system: 
Given the initial guess \(\psi_\Gamma^{0}\), compute for \(k=1,2,\dots\)
\begin{subequations}\label{eq:cont_cpl_laplace}
\begin{equation}
	\begin{aligned}
		c \frac{\partial\psi}{\partial t}  - K\frac{\partial^2}{\partial z^2} \psi &= 0, \\
		\hat{\psi}^k(s, -L) &= 0, \\
		\hat{\psi}^k(s, 0) &= \psi_\Gamma^{k-1},
	\end{aligned}
\end{equation}
\begin{equation}
	s\hat{h}^k = \hat{v}_s^k = -K \left(\partial_z \hat{\psi}^k(s,0) + 1\right),         
\end{equation}
\begin{equation}
	\hat{\psi}_\Gamma^k = \omega \hat{h}^k + (1-\omega) \hat{\psi}_\Gamma^{k-1}.
\end{equation}
\end{subequations}

Now we study the error of each iteration with respect to the true solution of the coupled problem, introducing
\begin{equation*}
	\hat{e}_1^k \coloneqq \hat{h}^{k} - \hat{h}, \quad \hat{e}_2^k \coloneqq \hat{\psi}^{k} - \hat{\psi}, \quad \hat{\gamma}^k \coloneqq \omega \hat{e}_1^k + (1-\omega) \hat{\gamma}^{k-1}.
\end{equation*}
The errors \(\hat{e}_1^k(s)\), \(\hat{e}_2^k(s,z)\) satisfy equation \eqref{eq:cont_cpl_laplace} with initial guess $\hat{e}_{\Gamma}^0$. We solve it by standard solution approaches for homogeneous linear PDEs (cf. \Cref{subsec:well_posedness}):
\begin{equation*}
	\hat{e}_1^k(s) = -\hat{\gamma}^{k-1} \sqrt{\frac{cK}{s}} \coth\left(\sqrt{\frac{cs}{K}} L\right), \ 
	\hat{e}_2^k(s,z) = \hat{\gamma}^{k-1} \frac{\sinh\left(\sqrt{\frac{cs}{K}} (z+L)\right)}{\sinh\left(\sqrt{\frac{cs}{K}}L\right)}.
\end{equation*}

This yields the convergence factor of the iterations in Laplace space
\begin{equation}\label{eq:continuous_rho}
    \rho(s,\omega) = \frac{\hat{e}_1^{k+1}}{\hat{e}_1^{k}} = \frac{\hat{e}_2^{k+1}}{\hat{e}_2^{k}} = \frac{\hat{\gamma}^k}{\hat{\gamma}^{k-1}} = 1 - \omega - \omega \sqrt{\frac{cK}{s}} \coth\left(\sqrt{\frac{cs}{K}} L\right).
\end{equation}
In particular, it holds that
\begin{equation*}
    \left\|\hat{e}_2^k(s,z)\right\|_2 = \left\|\rho^k \hat{e}_2^0(s,z)\right\|_2 \leq |\rho|^k \left\|\hat{e}_2^0(s,z)\right\|_2.
\end{equation*}
Parseval's theorem relates this result in Laplace space back to the time domain:
\begin{equation*}
    \left\|e^{-\eta t}e_2^k(t,z)\right\|_2 \leq |\rho|^k \left\|e^{-\eta t} e_2^0(t,z)\right\|_2.
\end{equation*}
Convergence in two iterations is obtained by requiring that \(\rho(s,\omega) \stackrel{!}{=} 0\), from which we obtain the optimal relaxation parameter
\begin{equation}
    \omega_\mathrm{opt}(s) = \frac{1}{1+\sqrt{\frac{cK}{s}}\coth\left(\sqrt{\frac{cs}{K}} L\right)}.
\end{equation}

The result indicates that the product of hydraulic capacity and conductivity \(cK\) is of importance for the convergence speed.
For \(\omega=1\), \(\mathrm{Re}(\rho) < 0\), consistent with continuous and discrete results for similar coupling iterations for the heat equation \parencite{gander.etal_2016,monge.birken_2018}.
Due to the dependence on \(s\), the application of these results to a discretization of the coupling iterations \eqref{eq:cont_cpl_iterations} and \eqref{eq:linear_continuous_cpl} is not straightforward.

\subsection{Discretization}

Analogously to \Cref{sec:discretization}, we now discretize \Cref{eq:linear_cont_richards} with linear finite elements in space and the implicit Euler method in time on an equidistant grid with grid size \(\dz=L/M\). This gives the following system of ODEs for the \(M-1\) interior nodal unknowns \(\psi_1(t),\dots, \psi_{M-1}(t)\):
\begin{equation*}
	\begin{aligned}
		c \sum_{j=1}^{M-1}\left\langle \phi_i, \phi_j \right\rangle \dot{\psi}_j^k 
	+ K \sum_{j=1}^{M-1}\left\langle \partial_z\phi_i, \partial_z\phi_j \right\rangle \psi_j^k &= f_{\Gamma,i}^{k-1} &\ \forall i=1,\dots,M-1, \\
		\psi_j(0) &= \psi_{0, j} &\ \forall j=1,\dots,M-1.
	\end{aligned}
\end{equation*}
The hat functions are denoted by \(\phi_i(z)\) and \( \psi_{0, j} \) are the coefficients that suitably approximate \( \psi_0(z) \) in \( \mathcal{V}^{M-1} \).
The right-hand side \(f_{\Gamma,i}^{k-1}\) contains the Dirichlet contribution from the interface boundary condition on \(\Gamma\), cf. \Cref{eq:richards_fgamma}.
In matrix-vector notation, we write this as 
\begin{equation*}\label{eq:semidiscr_richards}
	\begin{aligned}
		M_{II} \boldsymbol{\dot{\psi}}_I^k + A_{II} \boldsymbol{\psi}_I^k &= \boldsymbol{f}_\Gamma^{k-1}, \\
		\boldsymbol{\psi}_I^k(0) &= \boldsymbol{\psi}_0,
	\end{aligned}
\end{equation*}
where \(M_{II}\) and \(A_{II}\) are the mass and stiffness matrices.
The subscripts \(._I\) and \(._\Gamma\) here refer to the distinction between contributions related to the interior \(\Omega\) and interface \(\Gamma\), respectively.
This notation borrows from the theory of domain decomposition \parencite[e.g.,][]{toselli.widlund_2005} and is analogous to \textcite{monge.birken_2018}.

The first coupling condition \eqref{eq:ibc_dirichlet} states that the capillary head at the interface \(\Gamma\) is equal to the water height.
We refer to the associated unknown of the Galerkin approximation as \(\psi_\Gamma(t)=\psi_M(t)\approx\psi(t,0)\).
Since we are using a nodal basis, \(\psi_\Gamma(t)\equiv h(t)\) and we can define the full vector of unknowns
\begin{equation*}
	\boldsymbol{\psi}(t) \coloneqq \begin{pmatrix}
		\psi_{1} \\ \vdots \\ \psi_{M-1} \\ \psi_{M}
	\end{pmatrix} = \begin{pmatrix}
		\psi_{1} \\ \vdots \\ \psi_{M-1} \\ h
	\end{pmatrix} = \begin{pmatrix}
		\boldsymbol{\psi}_I \\ \psi_\Gamma
	\end{pmatrix} \in \mathbb{R}^{M}.
\end{equation*}

For the Dirichlet boundary condition, we obtain that

\begin{equation*}
	\boldsymbol{f}_\Gamma^k = - M_{I\Gamma} \dot{\psi}_\Gamma^{k} - A_{I\Gamma} \psi_\Gamma^{k},
\end{equation*}
with the vectors \(M_{I\Gamma},A_{I\Gamma}\in\mathbb{R}^{M-1}\) given by
\begin{equation}\label{eq:def_igamma}
	(M_{I\Gamma})_{j} \coloneqq \left\langle \phi_i, \phi_{M} \right\rangle, \quad (A_{I\Gamma})_j \coloneqq \left\langle \partial_z\phi_i, \partial_z\phi_{M} \right\rangle, \quad j=1,\dots,M-1.
\end{equation}

Discretizing in time with the implicit Euler method finally yields the fully discrete system of equations for linearized, one-dimensional subsurface flow:
\begin{equation}\label{eq:discrete_richards}
	\left(M_{II} + \Delta t A_{II}\right) \boldsymbol{\psi}_{I}^{n, k}  
	=-(M_{I\Gamma} + \Delta t A_{I\Gamma}) {\psi}_{\Gamma}^{n, k-1}
	+ M_{II} \boldsymbol{\psi}_{I}^{n-1} + M_{I\Gamma} \psi_\Gamma^{n-1}.
\end{equation}

For the second coupling condition \eqref{eq:ibc_neumann}, we aim to obtain an expression of \( v_s \) in terms of the weak form of Richards' equation, consistent with the FE discretization.
This can be done using Green's first identity\footnote{Note that this weak reformulation of the surface flux via Green's first identity is not possible in the nonlinear case.}, cf. \textcite{monge.birken_2018,toselli.widlund_2005}, giving:

\begin{equation*}
	\begin{aligned}
		v_s(t) &= 
	-\sum_{j=1}^{M} \dot{\psi}_j \underbrace{c\left\langle \phi_j,\phi_{M}\right\rangle}_{\eqqcolon {(M_{\Gamma})}_j} 
	- \sum_{j=1}^{M} \psi_j \underbrace{K\left\langle\partial_z\phi_j,\partial_z\phi_{M} \right\rangle}_{\eqqcolon {(A_{\Gamma})}_j} - K \underbrace{\left\langle 1,\partial_z \phi_{M} \right\rangle}_{=1}\\
	&= -M_{\Gamma}\boldsymbol{\dot{\psi}} -  A_{\Gamma} \boldsymbol{\psi} - K.
	\end{aligned}
\end{equation*}
We can split \(M_{\Gamma}\) and \(A_{\Gamma}\) into two parts
\begin{equation*}
	M_\Gamma =
	\begin{bmatrix}
		M_{\Gamma I} & M_{\Gamma \Gamma}
	\end{bmatrix}, \quad
	A_\Gamma =
	\begin{bmatrix}
		A_{\Gamma I} & A_{\Gamma \Gamma}
	\end{bmatrix},
\end{equation*}
where
\begin{equation}\label{eq:def_ggamma}
	\begin{aligned}
		M_{\Gamma I} &= {(M_{I \Gamma})}^T, \ &M_{\Gamma \Gamma}=c\left\langle \phi_M,\phi_{M}\right\rangle, \\
		A_{\Gamma I} &= {(A_{I \Gamma})}^T, \ &A_{\Gamma \Gamma}=K\left\langle\partial_z\phi_M,\partial_z\phi_{M} \right\rangle.
	\end{aligned}
\end{equation}
The semi-discrete equation for the water height \eqref{eq:linear_cont_swe} thus becomes
\begin{equation*}
	\dot{\tilde{h}}^k = \dot{\tilde{\psi}}_\Gamma^k = -\left(M_{\Gamma I} \boldsymbol{\dot{\psi}}^k_I + M_{\Gamma\Gamma} \dot{\psi}_\Gamma^{k-1} +  A_{\Gamma I} \boldsymbol{\psi}^k_I + A_{\Gamma\Gamma} \psi_\Gamma^{k-1} + K\right).
\end{equation*}

After discretization in time we are left with
\begin{equation}\label{eq:discrete_swe}
	\begin{split}
		\tilde{\psi}_\Gamma^{n,k} 
		= &- (M_{\Gamma I } + \Delta t A_{\Gamma I}) \boldsymbol{\psi}^{n, k}_I 
		- (M_{\Gamma\Gamma} + \Delta t A_{\Gamma\Gamma}) \psi_\Gamma^{n, k-1}  \\
		&+ M_{\Gamma I } \boldsymbol{\psi}^{n-1}_I + (1 + M_{\Gamma\Gamma}) \psi_\Gamma^{n-1} - \dt K
	\end{split}
\end{equation}
for the surface flow.
The resulting algorithm is given in \Cref{alg:dnwr_ideal_discrete}.

\begin{algorithm}
    \caption{Coupling iterations for the fully discrete model problem.}\label{alg:dnwr_ideal_discrete}
    \begin{algorithmic}
		\State $n\gets 1$, $\psi_\Gamma^{0}\gets h_0$
		\While{$n\leq N+1$}
		\State $k\gets 1$, ${\psi_\Gamma}^{n, 0}\gets {\psi_\Gamma}^{n-1}$
		\While{$\left|\tilde{\psi}_\Gamma^{n,k} - \psi_\Gamma^{n,k-1}\right|\geq \text{TOL}$}
			\State \textbf{Solve} \Cref{eq:discrete_richards} $\to$ $\boldsymbol{\psi}_I^{n,k}$ given $\psi_\Gamma^{n,k-1}$
			\State \textbf{Solve} \Cref{eq:discrete_swe} \(\to\) \(\tilde{\psi}_\Gamma^{n,k}\) given \(\boldsymbol{\psi}_I^{n,k}\)
			\State \textbf{Relaxation step:} \(\psi_\Gamma^{n,k} = \omega \tilde{\psi}_\Gamma^{n,k} + (1-\omega) \psi_\Gamma^{n,k-1}\)
			\State{\(k \gets k+1\)}
		\EndWhile{}
		\State{$n\gets n+1$}
		\State{$\boldsymbol{\psi}^{n-1}\gets\boldsymbol{\psi}^{n-1, k-1}$}
		\EndWhile{}
      \end{algorithmic}
\end{algorithm}

\subsection{Fully Discrete Analysis}\label{subsec:optimal_omega}

We now want to find the convergence rate and optimal relaxation parameter of the fully discrete \Cref{alg:dnwr_ideal_discrete}.
As in \textcite{monge.birken_2018}, we achieve this by expressing \( \psi_\Gamma^{n, k} \) in terms of \( \psi_\Gamma^{n, k-1} \), \( \omega \), and values from the previous time step \( t^{n-1} \).
One first solves \Cref{eq:discrete_richards} for \(\boldsymbol{\psi}_I^{n,k}\) and plugs the result into \Cref{eq:discrete_swe}:
\begin{equation}\label{eq:update_height_norelax}
	\tilde{\psi}_\Gamma^{n,k} = S \psi_\Gamma^{n,k-1} + \xi^{n-1}. 
\end{equation}
Here \(\xi^{n-1}\) summarizes terms which only depend on the previous time step \(t^{n-1}\), while \(S\) is defined as\footnote{As also observed in \textcite{monge.birken_2018}, \( S \) takes the form of a Schur complement.}
\begin{equation*}
	S \coloneqq 
		(M_{\Gamma I } + \dt A_{\Gamma I}) \left(M_{II} + \dt A_{II}\right)^{-1}(M_{I\Gamma} + \dt A_{I\Gamma})
		- (M_{\Gamma\Gamma} + \dt A_{\Gamma\Gamma}).
\end{equation*}
Applying the relaxation step, we obtain
\begin{equation}\label{eq:update_height_relax}
	\psi_\Gamma^{n, k} = \underbrace{\left(\omega S + (1-\omega)\right)}_{\eqqcolon \Sigma(\omega)} \psi_\Gamma^{n, k-1} + \omega \xi^{n-1}. 
\end{equation}

In this model problem, the iteration matrix \( \Sigma(\omega) \) is a scalar, which reduces the optimality condition to \( \Sigma(\omega_\mathrm{opt}) \stackrel{!}{=} 0 \).
Therefore,
\begin{equation*}
	\omega_\mathrm{opt} = \frac{1}{1-S}.
\end{equation*}

Computing \(S\) involves the inverse of \( M_{II} + \dt A_{II} \).
Linear finite elements with a nodal basis, along with our assumed boundary conditions, yield that this is a symmetric, tridiagonal Toeplitz matrix given by
\begin{equation*}
	\begin{gathered}
		M_{II} + \dt A_{II} = \begin{bmatrix}
			a & b \\
			b & a & b \\
			& \ddots & \ddots & \ddots \\
			& & b & a & b \\
			& & & b & a
		\end{bmatrix} \in \mathbb{R}^{(M-1)\times (M-1)},\\
		a \coloneqq \frac{2}{3} c\dz + \frac{2K\dt}{\dz}, \\
		b \coloneqq \frac{1}{6} c\dz - \frac{K\dt}{\dz}.
	\end{gathered}
\end{equation*}
The sparsity of this matrix is not preserved in its inverse, but we can follow \textcite{monge.birken_2018} to mitigate this problem.
Since the matrix is Toeplitz,
\begin{equation*}
	{(M_{II} + \dt A_{II})}^{-1} = {(V \Lambda V^T)}^{-1} = V \Lambda^{-1} V^T.
\end{equation*}
Here, \( \Lambda=\mathrm{diag}(\boldsymbol{\lambda}) \), \( \boldsymbol{\lambda} = (\lambda_1, \dots, \lambda_{M-1})^T \)  is a diagonal matrix containing the eigenvalues of \( M_{II} + \dt A_{II} \) and \( V \) is an orthogonal matrix with the right eigenvectors of \( M_{II} + \dt A_{II} \).
Both are known for Toeplitz matrices:
\begin{equation}\label{eq:toeplitz_eigen}
	\begin{gathered}
		\lambda_j = a - 2b \cos\left(\frac{j\pi}{M}\right) \\
		V_{ij} = \frac{1}{\sqrt{\sum_{k=1}^{M-1} \sin^2\left(\frac{k\pi}{M}\right)}} \sin\left(\frac{ij\pi}{M}\right).
	\end{gathered}
\end{equation}

Evaluating the integrals in \Cref{eq:def_igamma}, we furthermore obtain that:
\begin{equation*}
	M_{I\Gamma} + \dt A_{I\Gamma} = {(0, \dots, 0, b)}^T.
\end{equation*}
Thus,
\begin{equation*}
	(M_{I \Gamma} + \dt A_{I \Gamma})^T (M_{II} + \dt A_{II})^{-1}(M_{I\Gamma} + \dt A_{I\Gamma}) = b^2 \alpha,
\end{equation*}
with
\begin{equation*}
	\alpha \coloneqq {\left({(M_{II} + \dt A_{II})}^{-1}\right)}_{M-1,M-1} = \sum_{j=1}^{M-1} \lambda_j^{-1} {V_{M-1,j}}^2.
\end{equation*}
Substituting \Cref{eq:toeplitz_eigen} gives
\begin{equation*}
	\alpha = \frac{\sum_{j=1}^{M-1} \sin^2\left(\frac{(M-1)j\pi}{M}\right) \left(a - 2b \cos\left(\frac{j\pi}{M}\right)\right)^{-1}}{\sum_{j=1}^{M} \sin^2\left(\frac{j\pi}{M}\right)}.
\end{equation*}

We can simplify this expression using the identities
\begin{equation*}
	\begin{aligned}
		\sin^2\left(\frac{(M-1)j\pi}{M}\right) &= \sin^2\left(\frac{j\pi}{M}\right) \quad \forall j\in\mathbb{N}, \\
		\sum_{j=1}^{M-1} \sin^2\left(\frac{j\pi}{M}\right) &= \frac{M}{2}.
	\end{aligned}
\end{equation*} 
Remembering that \( \dz = L/M \), we thus obtain:
\begin{equation}
	\alpha = \frac{\dz}{L} \sum_{j=1}^{M-1} \frac{\sin^2\left(\frac{j\pi\dz}{L} \right)}{\frac{a}{2} - b \cos\left(\frac{j\pi\dz}{L}\right)}.
\end{equation}

Finally, we note that we can simplify the subtrahend by evaluating the integrals in \Cref{eq:def_ggamma}:
\begin{equation*}
	M_{\Gamma \Gamma} + \dt A_{\Gamma\Gamma} = \frac{a}{2}.
\end{equation*}
With this we get that
\begin{equation}\label{eq:discr_analysis_result}
	S = b^2 \alpha - \frac{a}{2}, \quad \omega_\mathrm{opt} = \frac{1}{1 + \frac{a}{2} - b^2 \alpha}.
\end{equation}

The first term in \(S\) is positive\footnote{Recall that \( (M_{\Gamma I } + \dt A_{\Gamma I}) = (M_{I\Gamma} + \dt A_{I\Gamma})^T \). If we define \( \boldsymbol{x} \coloneqq (M_{I\Gamma} + \dt A_{I\Gamma}) \neq 0 \), \( A \coloneqq \left( M_{II} + \dt A_{II}\right)^{-1} \), the first term in \( S \) corresponds to \( \|\boldsymbol{x}\|_A^2 > 0 \).}, and thus in particular \(\alpha>0\).
One result of the continuous analysis was that the real part of the convergence factor \(\mathrm{Re}(\rho)\) is negative, which indicates that \(S<0\) should hold in the discrete analysis.
However, we are unable to show that this holds in general.

\section{Numerical Results}\label{sec:numerical_results}

\subsection{Implementation and Setup}

To verify and test our analysis, we have implemented both the nonlinear 2D-1D coupling iteration \eqref{eq:cont_cpl_iterations} and the linear 1D-0D coupling iteration \eqref{eq:linear_continuous_cpl}.
Our code makes use of the DUNE Numerics Toolbox
\parencite{bastian.etal_2021,dedner.etal_2020} for the space and time discretization.
DUNE is a general purpose open source framework for solving PDEs with grid
based numerical schemes. It is written in C++ but
provides Python bindings for ease of use, which are the basis of our simulations.
The models are defined using the unified form language (UFL) \parencite{alnaes.etal_2014} which
simplifies the definition of weak formulations of PDEs via a domain specific language.

The solvers are coupled with preCICE \parencite{chourdakis.etal_2022}, an open-source C++ coupling library which also provides Python bindings.
preCICE natively supports sequential coupling (via the \emph{serial-implicit} coupling scheme) with constant acceleration by a relaxation parameter, as required for the coupling iterations.
Our implementation is available as open source code.\footnote{Code available at: \url{https://gitlab.maths.lu.se/nuan/projects/coupled-hydrology/-/releases/1.0.0}.}

We use the termination criterion \eqref{eq:conv_crit} with \(\mathrm{TOL} = \num{e-8}\) in all experiments.
We also define an experimental convergence factor

\begin{equation}\label{eq:conv_estimate}
    \expcvg_n = \frac{1}{K_n-2}\sum_{k=2}^{K_n-1}\frac{\mathrm{res}^{n,k}}{\mathrm{res}^{n,k-1}}\approx\rho,
\end{equation}
where \(K_n\) is the amount of iterations in time step \(n\) until \Cref{eq:conv_crit} is met.
In case of sequential coupling in the 1D-0D model problem, it holds that:

\begin{equation}\label{eq:cr_motivation}
    \expcvg_n = |\Sigma(\omega)|.
\end{equation}

\noindent
See \Cref{app:convergence_estimate} for a derivation of this result.
We will also make use of the time-averaged experimental convergence rate

\begin{equation*}
    \expcvg=\frac{1}{N}\sum_{n=1}^{N} \expcvg_n.
\end{equation*}

\subsection{Verification of Linear 1D-0D Analysis}

In this section we verify and study results from the discrete analysis presented in \Cref{sec:linear_analysis}.
The theoretical analysis was carried out for a single time step.
We mirror this by choosing the final simulation time \(T=\dt\), i.e., \(\expcvg \equiv \expcvg_1\).
All simulations in this section use the initial conditions \(\psi_0(z)=1-z/L\), \(h_0=0\).

\begin{figure}
    \centering
    \begin{subfigure}{0.31\textwidth}
        \includegraphics[height=9em]{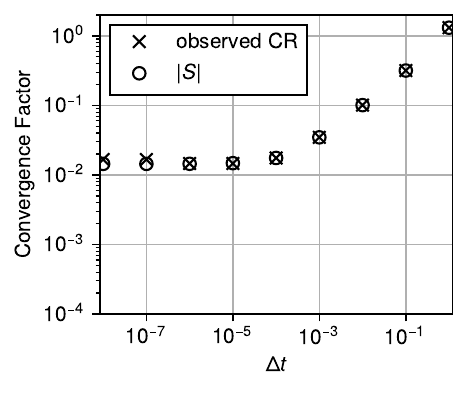}
        \caption{\(\dz=1/20\), \(\omega=1\)}
        \label{fig:cvg_no_relaxation_coarse}
    \end{subfigure}
    \hfill
    \begin{subfigure}{0.31\textwidth}
        \includegraphics[height=9em]{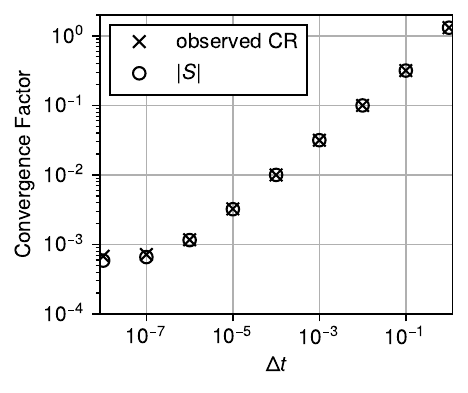}
        \caption{\(\dz=1/500\), \(\omega=1\)}
        \label{fig:cvg_no_relaxation_fine}
    \end{subfigure}
    \hfill
    \begin{subfigure}{0.31\textwidth}
        \includegraphics[height=9em]{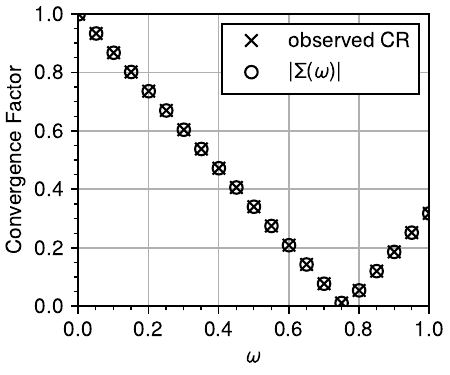}
        \caption{\(\dz=1/20\), varying \(\omega\).}
        \label{fig:cvg_relaxation}
    \end{subfigure}
    \hfill
    \caption{Verification of fully discrete analysis: Comparison of the theoretical and experimental convergence factors for varying grid and time step sizes. \label{fig:analysis_verification}}
\end{figure}

We first verify the theoretical convergence factor in case the relaxation parameter \(\omega=1\), i.e., where \(\Sigma \equiv S\).
For this test, we set \(c=K=L=1\) and varied the discretization parameters \(\dt,\dz\).
The results are displayed in \Cref{fig:cvg_no_relaxation_coarse,fig:cvg_no_relaxation_fine} for two different values of \(\dz\).
The difference between the experimentally observed and theoretical convergence factor is on the order of \(10^{-11}\) to \(10^{-9}\).
Thus, the analysis agrees well with observed convergence rates.

Next we check that the optimal relaxation parameter \(\omega_\mathrm{opt}=1/(1-S)\) is correctly predicted by our analysis.
We have verified this for a range of values and include figures here for the case \(c=K=L=1\) and \(\dt=\dz=10^{-1}\).
\Cref{fig:cvg_relaxation} shows that the results are in agreement for \(\omega\in[0,1]\).
The linear dependence of \(\Sigma\) on \(\omega\) is clearly visible.
In particular, choosing \(\omega=\omega_\mathrm{opt}\) gives a very small convergence factor (and convergence in two iterations with \(\mathrm{TOL}=\num{e-8}\)).
These experiments confirm that \Cref{eq:discr_analysis_result} captures the convergence properties of the linearized, fully discrete coupling iteration (\Cref{alg:dnwr_ideal_discrete}).

\begin{figure}
    \centering
    \begin{subfigure}{0.3\textwidth}
        \centering
        \includegraphics[height=11em]{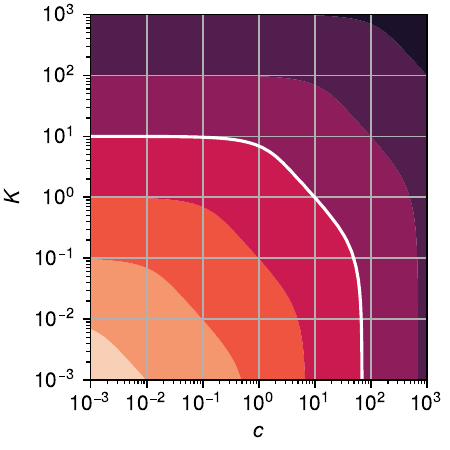}
        \caption{Varying physics, fixed \(\dz=1/20, \dt=\num{e-1}\).}
        \label{fig:factor_varying_physics_coarse}
    \end{subfigure}
    \hfill
    \begin{subfigure}{0.3\textwidth}
        \centering
        \includegraphics[height=11em]{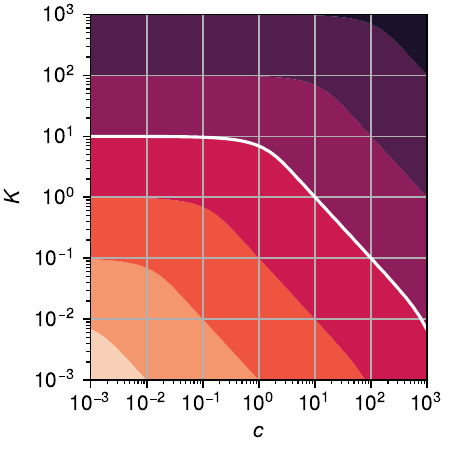}
        \caption{Varying physics, fixed \(\dz=1/500, \dt=\num{e-1}\).}
        \label{fig:factor_varying_physics_fine}
    \end{subfigure}
    \hfill
    \begin{subfigure}{0.31\textwidth}
        \centering
        \includegraphics[height=11em]{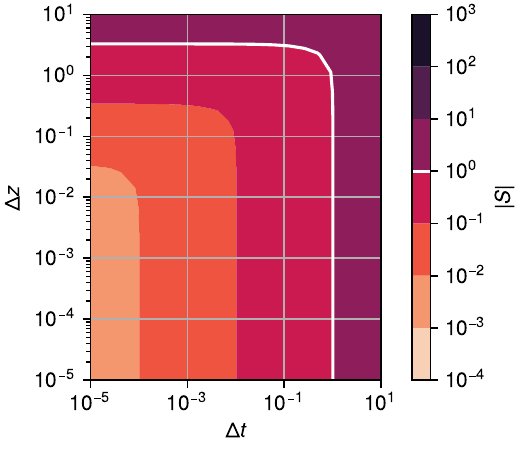}
        \caption{Fixed \(c=K=1\), varying grid sizes.}
        \label{fig:factor_varying_grids}
    \end{subfigure}
    \hfill
    \caption{Convergence factor \(|S|\) for varying hydraulic conductivity $K$, hydraulic capacity $c$, and grid sizes. 
    The white contour line for \(|S|=1\) marks the border between converging and diverging coupling iterations. \label{fig:factor_analysis}}
\end{figure}

Now that we have verified the results of our analysis, we can use \Cref{eq:discr_analysis_result} to determine the convergence properties of the sequential coupling iteration. \Cref{fig:factor_analysis} shows results for the convergence factor \(|S|\), for a wide range of material parameters \(c, K\) and discretization parameters \(\dt, \dz\), respectively. 
As described in \Cref{sec:linear_analysis}, we expect \(S\) to be negative.
Indeed, \(S<0\) for all parameter choices we tested.
The white contour lines mark the transition from convergent to divergent coupling iterations, i.e., \(|S|=1\). The values of $|S|$ vary from \numrange{e-4}{e3}.

\begin{figure}
    \centering
    \begin{subfigure}{0.32\textwidth}
        \includegraphics[height=11em]{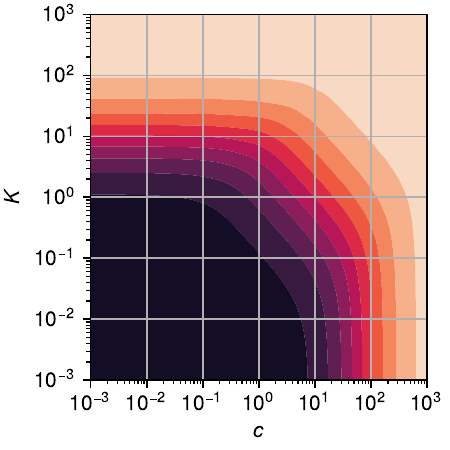}
        \caption{\(\dz=1/20, \dt=\num{e-1}\)}
        \label{fig:varying_physics_coarse}
    \end{subfigure}
    \hfill
    \begin{subfigure}{0.32\textwidth}
        \includegraphics[height=11em]{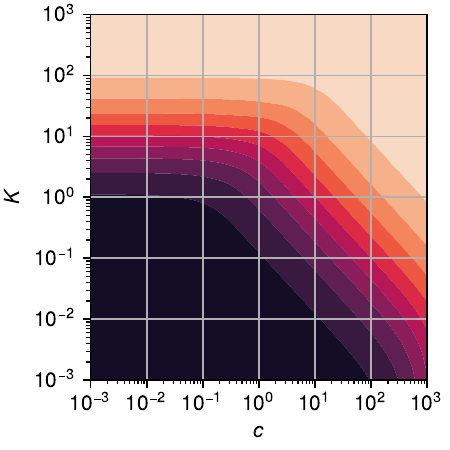}
        \caption{\(\dz=1/500, \dt=\num{e-1}\)}
        \label{fig:varying_physics_fine}
    \end{subfigure}
    \hfill
    \begin{subfigure}{0.32\textwidth}
        \includegraphics[height=11em]{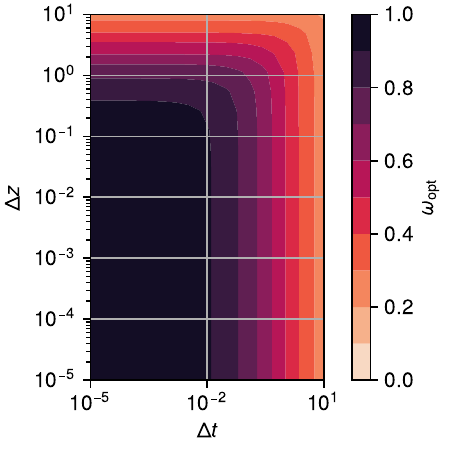}
        \caption{\(c=K=1\)}
        \label{fig:varying_grids}
    \end{subfigure}
    \hfill
    \caption{\(\omega_\mathrm{opt}\) for varying hydraulic conductivity $K$, hydraulic capacity $c$, and grid sizes. \label{fig:omega_opt_analysis}}
\end{figure}

\Cref{fig:factor_varying_physics_coarse,fig:factor_varying_physics_fine} show the convergence factor for fixed grid sizes but varying physical constants \(c\) and \(K\).
Since these constants can take on various orders of magnitude due to the nonlinear constitutive relations \eqref{eq:genuchten}, we vary both \(c\) and \(K\) between \num{e-3} and \num{e3}.
\Cref{fig:factor_varying_physics_coarse} plots the results for a coarser grid \(\dz=1/20\), whereas \(\dz=1/500\) in \Cref{fig:varying_physics_fine}.
In both cases, \(\dt=1/10\).

One can distinguish three regimes:
In the case where \(c \ll K\), the convergence factor scales with \(K\); this can be seen in both \Cref{fig:factor_varying_physics_coarse,fig:factor_varying_physics_fine}. Vice versa, \(c \gg K\) yields a convergence factor that only depends on \(c\), as can be seen in the lower right corner of \Cref{fig:varying_physics_coarse}.
When the hydraulic conductivity and capacity have similar orders of magnitude, the diagonal slope of the contour lines shows a dependence on the product \(c K\).
Decreasing the grid size \(\dz\) shifts and scales the results, increasing the size of this last regime.
Note that this relationship between the convergence factor and \(c K\) is especially clear in the result from the continuous analysis \eqref{eq:continuous_rho}:
If the product of hydraulic conductivity and capacity increases while their ratio stays the same, \(|\rho(s,\omega=1)|\) will also increase.
Since a numerical discretization should approach the continuous solution in the limit of small grid sizes, it is reasonable that a smaller grid size \(\dz\), as in \Cref{fig:cvg_no_relaxation_fine}, increases the size of the regime where this behavior is observed.

\Cref{fig:factor_varying_grids} displays the convergence factor for \(c=K=1\) with varying discretization parameters \(\dt,\dz\).
Again we see three regimes:
For \(\dt\ll\dz\), \(S\) mainly depends on \(\dz\), representing a situation where the numerical error in space dominates.
Further decreasing the time step size will thus not help improve convergence of the iterations.
Vice versa, for \(\dt\gg\dz\), the contour lines show that the most effective way of decreasing \(|S|\) is by improving the resolution in time.
The transition between these two extremes follows \(\dz = \sqrt{\dt}\), which is expected for the heat equation.

These figures show that it is possible to obtain non-convergent coupling iterations for coarse resolutions in time and space, combined with large values for the hydraulic capacity and conductivity.

\Cref{fig:omega_opt_analysis} illustrates parameter choices where the coupling iteration can be accelerated by adjusting the relaxation parameter.
\Cref{fig:varying_physics_coarse,fig:varying_physics_fine} show the dependence of \(\omega_\mathrm{opt}\) on \(c\) and \(K\) for fixed \(\dt=1/10\) and \(\dz=1/20\) and \(\dz=1/500\), respectively.
The dependence on discretization parameters for fixed \(c=K=1\) is displayed in \Cref{fig:varying_grids}.

The plots in \Cref{fig:omega_opt_analysis} exhibit the same regimes as described for the convergence factor above.
They illustrate that \(\omega_\mathrm{opt}\) can take on values from close to 0 all the way up to 1.
Relaxation ceases to be useful when \(|S|\) is very small or very large, corresponding to the dark and light regions in the figure (where \(\omega_\mathrm{opt}\approx 1\) and \(\omega_\mathrm{opt}\approx 0\)).
In the former case, the iterations terminate quickly and one can focus on solver, instead of coupler, development.
In the latter case, the coupling problem itself should be adjusted to make the iterations more robust (e.g., by adjusting interface boundary conditions).

\subsection{Results with Nonlinear 2D-1D Code}

We now turn to investigate how well the linear analysis is able to explain
the behavior in simulations with the nonlinear models presented in \Cref{sec:governing_equations,sec:discretization}.
Past numerical studies for coupled surface-subsurface modeling relied on comparing to observations or intermodel comparisons with benchmark problems \parencite{maxwell.etal_2014,sulis.etal_2010,kollet.maxwell_2006,vanderkwaak_1999}.
These all start off with zero water height in the surface solver.
Since we are interested in the case where the shallow water solver is active, we have not implemented a boundary condition switching procedure and instead modified existing examples to fit our setup.

We look at two different test cases based on the benchmark problems from \textcite{schneid_2000,maxwell.etal_2014}, which differ in spatial scale and make use of realistic materials.
All van Genuchten parameters used in the experiments are given in \Cref{tab:van_genuchten_params}.

The hydraulic conductivity \(K\) has its maximum at \(K_s\) for \(\psi\geq 0\).
One can use the van Genuchten equations \eqref{eq:genuchten},\eqref{eq:capacity_definition} to additionally determine theoretical maxima of the capacity \(c\).
We have reported these in the last two rows of \Cref{tab:van_genuchten_params}.
The values obtained for these materials are small, indicating that one can expect fast convergence of the iterations according to the linear analysis (cf. \Cref{fig:factor_varying_physics_coarse,fig:factor_varying_physics_fine}).

The convergence behavior of the coupling iterations is measured using preCICE output data and the experimental convergence factor \(\expcvg_n\).
We compare this with an expected convergence factor \(|S|\) for each model output time step: we compute the spatial means of \(K(\psi^n)\) and \(c(\psi^n)\), \(\bar{K}^n, \bar{c}^n\), and plug them into \Cref{eq:discr_analysis_result}.

It is reasonable to assume that convergence of the coupling iterations is driven by the material properties in the vicinity of the interface, which would motivate estimating \(|S|\) only on, or near, \(\Gamma\).
However, \(c=0\) where \(\psi\geq 0\) and by the coupling condition \eqref{eq:ibc_dirichlet} \(\psi\geq 0\) on \(\Gamma\), rendering the linear model invalid.
For this reason we stick with taking an average over the whole domain, thus overestimating the influence of deeper soil regions.

\begin{table}[]
    \resizebox{\textwidth}{!}{%
    \begin{tabular}{|ll|l|l|l|l|}
    \hline
    \textbf{Parameter} & \textbf{} & \textbf{Unit} & \textbf{Beit-Netofa Clay} & \textbf{Silt Loam} & \textbf{Sandy Loam} \\ \hline
    Alpha & $\alpha$ & \unit{m^{-1}} & \num{0.152} & \num{0.423} & \num{100.0} \\
    Pore-size distributions & $n_G$ & - & \num{1.17} & \num{2.06} & \num{2.0} \\
    Residual water content & $\theta_R$ & - & \num{0.0} & \num{0.131} & \num{0.2} \\
    Saturated water content & $\theta_S$ & - & \num{0.446} & \num{0.396} & \num{1.0} \\
    Sat. hydraulic conductivity & $K_s$ & \unit{m \per s} & \num{9.49e-9} & \num{5.74e-7} & \(1.16 \times 10^{\{-5,-6,-7\}}\) \\
    Max. hydraulic capacity & $c$ & \unit{m^{-1}} & \num{7.45e-3} & \num{0.0449} & \num{30.6} \\ \hline
    \end{tabular}
    }
    \caption{The van Genuchten parameters from \textcite{list.radu_2016,maxwell.etal_2014} used in the nonlinear numerical experiments. The last two rows state maxima for \(K\) and \(c\) for all \(\psi\) according to the van Genuchten model.}
    \label{tab:van_genuchten_params}
\end{table}

\subsubsection{Small-Scale Experiment: Drainage Trench}

This test case is based on a drainage trench example for subsurface modeling \parencite{list.radu_2016,schneid_2000}.
The original problem models the vertical cross-section of a groundwater domain that recharges from a drainage trench covering a part of the surface.
\textcite{list.radu_2016} study this case for silt loam and Beit-Netofa clay.

We look at a variation of this test by covering the whole ground with a flat water surface which fills up with the same inflow rate as in the original setup:
\begin{equation*}
    r(t) = \begin{cases}
        \qty{10}{cm.h^{-1}} & t \leq \qty{2}{h},\\
        0 & t > \qty{2}{h}.
    \end{cases}
\end{equation*}
The domain is two meters wide and flat, which is why we used the small scale shallow water model \eqref{eq:swe} here.

We consider three soil variants for the subsurface domain:
homogeneous silt loam, homogeneous Beit-Netofa clay, and a heterogeneous case.
In the latter, the left half of the domain is made up of silt loam, whereas the right half is made up of Beit-Netofa clay.
We smoothen the transition between the two soils using an interpolation with a scaled \(\tanh\):
\begin{equation*}
    \beta(x) = \frac{\tanh\left(4 \left(x - L_x / 2\right)\right) + 1}{2}.
\end{equation*}
We interpolate all parameters in the van Genuchten--Mualem model \eqref{eq:genuchten} with this approach.
For instance, the saturated hydraulic conductivity becomes:
\begin{equation*}
    K_s(x) = (1 - \beta(x)) K_{s,\mathrm{loam}} + \beta(x) K_{s,\mathrm{clay}}.
\end{equation*}

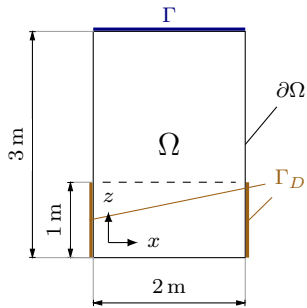
\begin{figure}
    \centering
    \begin{tikzpicture}[circuit ee IEC]
    \tikzstyle{doublefullarrow}=[<->, >=triangle 60, very thin]
    \tikzset{>=latex}
    \draw[draw=black] (0.0, 0.0) rectangle ++(2.0, 3.0);
    
    \node[black] at (1.0, 1.5) {\large $\Omega$};
    \node[anchor=south] (dOmega) at (2.6, 2) {\footnotesize $\partial\Omega$};
    \draw[black, -] (2.0, 1.5) -- (dOmega);
    
    \node[groundColor, anchor=south] (GammaD) at (2.6, 0.8) {\footnotesize $\Gamma_D$};
    \draw[groundColor, -] (2.05, 0.5) -- (GammaD);
    \draw[groundColor, -] (-0.05, 0.5) -- (GammaD);

    \node[] (A) at (0, 0) {};
    \node[] (B) at (0, 3) {};
    \node[] (C) at (2, 3) {};
    \node[] (D) at (2, 0) {};
    \node[] (E) at (0, 1) {};
    \node[] (F) at (2, 1) {};

    \draw[groundColor, very thick] ($(A) + (-0.03, 0)$) -- ($(E) + (-0.03, 0)$);
    \draw[groundColor, very thick] ($(D) + (0.03, 0)$) -- ($(F) + (0.03, 0)$);

    \draw[waterColor, very thick] ($(B) + (0, 0.03)$) -- ($(C) + (0, 0.03)$);
    \node[waterColor, anchor=south] (Gamma) at (1, 3) {\footnotesize $\Gamma$};

    \node[shape=coordinate] (CSorigin) at ($(A) + (0.2, 0.2)$) {};
    \node[] (CSz) at ($(CSorigin) + (0, 0.6)$) {\footnotesize $z$};
    \node[] (CSx) at ($(CSorigin) + (0.6, 0)$) {\footnotesize $x$};
    \draw[black, ->] (CSorigin) -- (CSz);
    \draw[black, ->] (CSorigin) -- (CSx);

    \dimline[extension start length=0.3, extension end length=0.3, label style={above=0.02ex,},] {($(A) + (-0.3, 0)$)} {($(E) + (-0.3, 0)$)} {\footnotesize \(\qty{1}{m}\)};
    \dimline[extension start length=0.27, extension end length=0.27, label style={above=0.02ex,},] {($(A) + (-0.8, 0)$)} {($(B) + (-0.8, 0)$)} {\footnotesize \(\qty{3}{m}\)};
    \dimline[extension start length=-0.3, extension end length=-0.3, label style={above=0.02ex,},] {($(A) + (0, -0.6)$)} {($(D) + (0, -0.6)$)} {\footnotesize \(\qty{2}{m}\)};
    \draw [dashed] (E) -- (F);
    \draw[ground] (1, 0.5) {};
\end{tikzpicture}
    \caption{Geometry for the drainage trench experiments. The dashed line marks the initial groundwater height. \label{fig:schneid_domain}}
\end{figure}

The geometry for the subsurface domain is given in \Cref{fig:schneid_domain}.
As initial and boundary conditions, we choose
\begin{equation*}
    \begin{aligned}
        \psi(0,\boldsymbol{x}) &= 1-z \quad\text{on } \Omega,\\
        \psi(t,\boldsymbol{x}) &= 1-z \quad\text{on } \Gamma_D,\\
        \boldsymbol{v} \cdot \boldsymbol{n} &= 0   \quad\text{on } \partial\Omega \setminus \{\Gamma_D \cup \Gamma\},\\
        h(0, x) &= \qty{1e-6}{m},\\
        \partial_x h(t, 0) &= \partial_x h(t, 2) = 0.
    \end{aligned}
\end{equation*}
This corresponds to no-flux boundary conditions on \(\Omega\) everywhere except for the Dirichlet boundary \(\Gamma_D\), and outflow boundary conditions on \(\Gamma\).
Since \Cref{eq:swe} is not defined for \(h=0\), we initialize the water height with a small, nonzero value.
In all experiments \(T=\qty{3}{h}\), \(\omega=1\), \(\dx=\qty{0.4}{m}\), \(\dz=\qty{0.375}{m}\), and \(\dt=\qty{36}{s}\).

\begin{figure}
    \centering
    \includegraphics[width=0.7\textwidth]{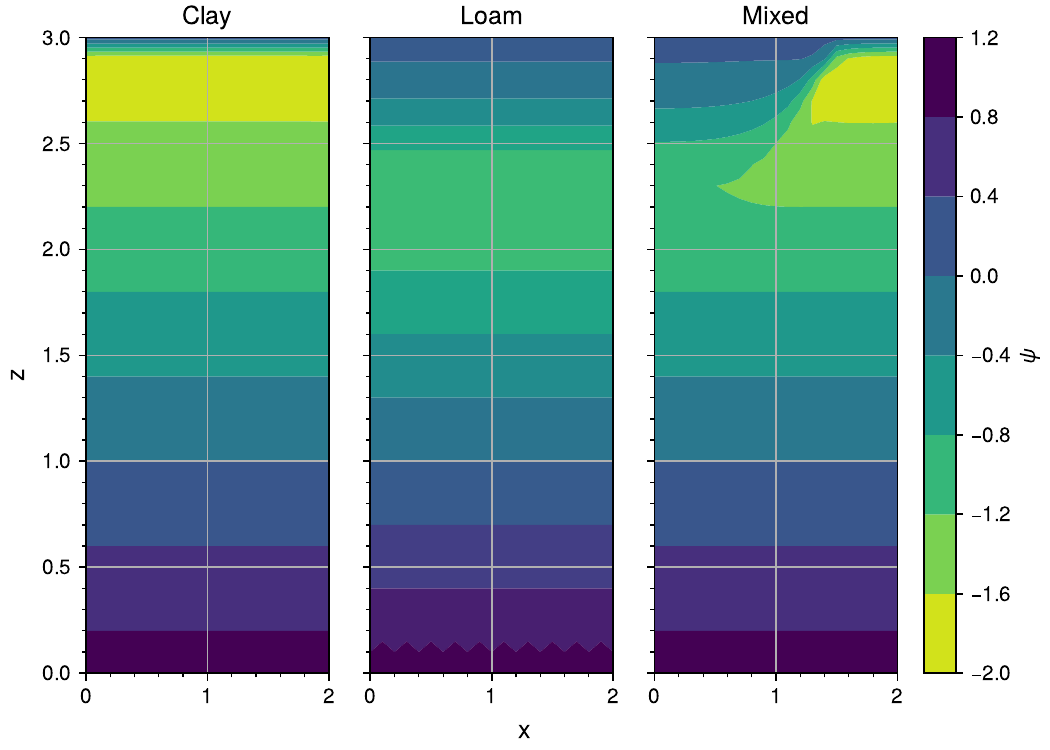}
    \caption{Capillary head $\psi(t, \boldsymbol{x})$ observed at the end of the drainage trench simulation ($t=\qty{3}{h}$). \label{fig:drainage_psi}}
\end{figure}

\Cref{fig:drainage_psi} shows simulation results for the capillary head \(\psi(t=T,\boldsymbol{x})\) at the end of the simulation for the three soil configurations.
Positive values \(\psi>0\) correspond to saturated soil.
The first two panels clearly show that clay takes up surface water more slowly than loam, in line with the results of \textcite{list.radu_2016}.
The third panel illustrates the nonlinearity of the model, yielding a complex saturation profile in the heterogeneous soil configuration.

\begin{figure}
    \centering
    \includegraphics[width=0.9\textwidth]{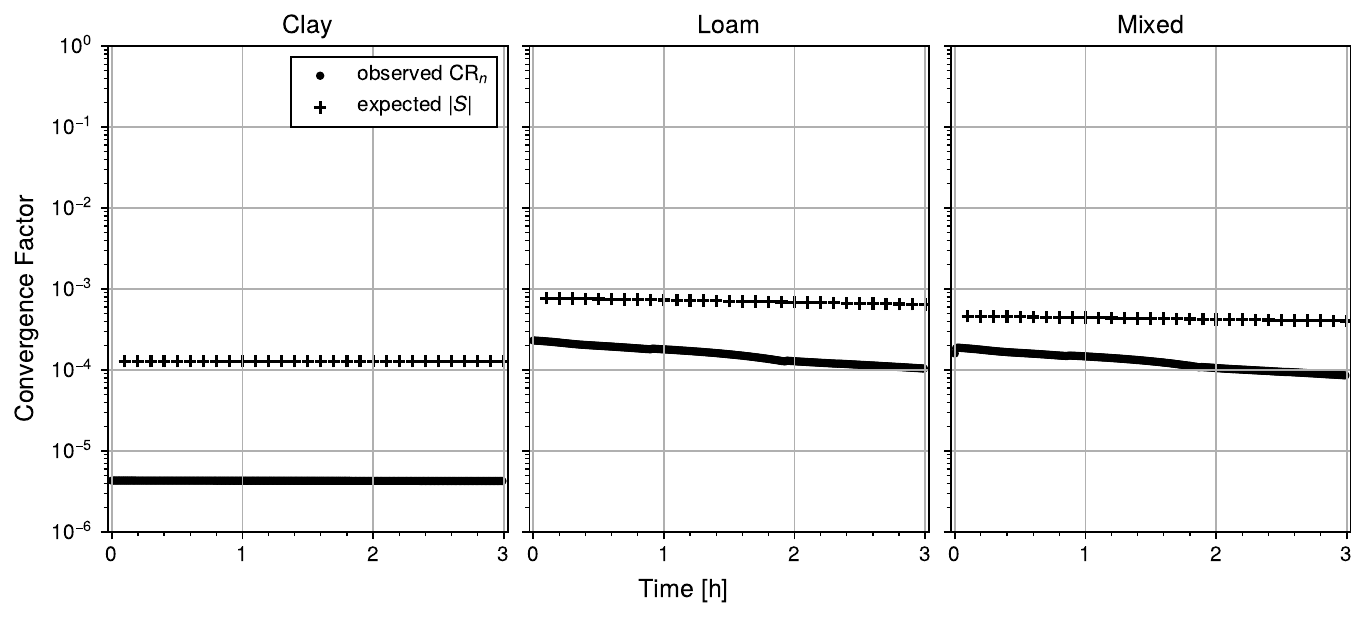}
    \caption{Observed vs. expected convergence factors for the drainage trench example. \label{fig:drainage_convergence_factors}}
\end{figure}

We now compare the observed convergence factors \(\expcvg_n\) to the expected convergence factor \(|S|\) over the course of a simulation.
\Cref{fig:drainage_convergence_factors} shows the results for the three soil configurations.
Note that we compute \(\expcvg_n\) in every time step while \(|S|\) is only computed every time we write an output file, in this case every ten time steps.

All observed convergence rates are very small, on the order of \numrange{e-6}{e-3}.
They are lowest for clay and highest for loam, while the heterogeneous soil is slightly below the loam results.
Since we use the same time step size and grid resolution for all three cases, this difference must be due to the material properties.
The convergence rate for Beit-Netofa clay stays approximately constant over time, while it visibly decreases for the other two cases.

The expected convergence factor correctly predicts the low convergence rates in all cases, along with the qualitative difference between soil types, and the time dependence.
The low convergence rate matches the small observed material parameters observed during the simulations, with the maximal values being \(c=\num{0.045}\), \(K=\num{5.74e-7}\) (both in the silt loam case).
Comparing with \Cref{fig:factor_analysis}, such small values correspond to the regime where \(|S|\ll 1\).
However, we see that the expected convergence factor is consistently above the observed values, particularly in case of clay.
In practice, this does not have a significant effect in terms of accelerating convergence:
The largest predicted convergence factor is around \num{e-3}, giving \(\omega_\mathrm{opt} > 0.999\) for all three experiments.
This means that relaxation is not necessary for fast convergence of the coupling iterations.

We can explain the time dependence using the van Genuchten--Mualem model and the discrete analysis results from \Cref{fig:factor_analysis}:
Over time, more parts of the soil are saturated, giving \(\psi>0\) and thus \(\lim_{n\to\infty} \bar{c}^n = 0\) while \(\lim_{n\to\infty} \bar{K}^n = K_s\).
Comparing this with \Cref{fig:factor_varying_physics_coarse,fig:factor_varying_physics_fine}, we see that smaller values of \(c\) with finite \(K\) indeed yield a smaller convergence factor.

A similar reasoning can explain why the predicted values are consistently higher than the observed values, and why this difference is largest for Beit-Netofa clay:
As seen in \Cref{fig:drainage_psi}, the soil has higher degrees of saturation close to the interface, speeding up convergence.
This behavior is not as prominent in our computation of \(|S|\) since we overestimate the influence of deeper, less saturated, soil regions (by taking the spatial means of \(c,K\)).
The saturated near-interface region is smallest for Beit-Netofa clay, amplifying the effect.

\begin{figure}
    \centering
    \includegraphics[width=0.9\textwidth]{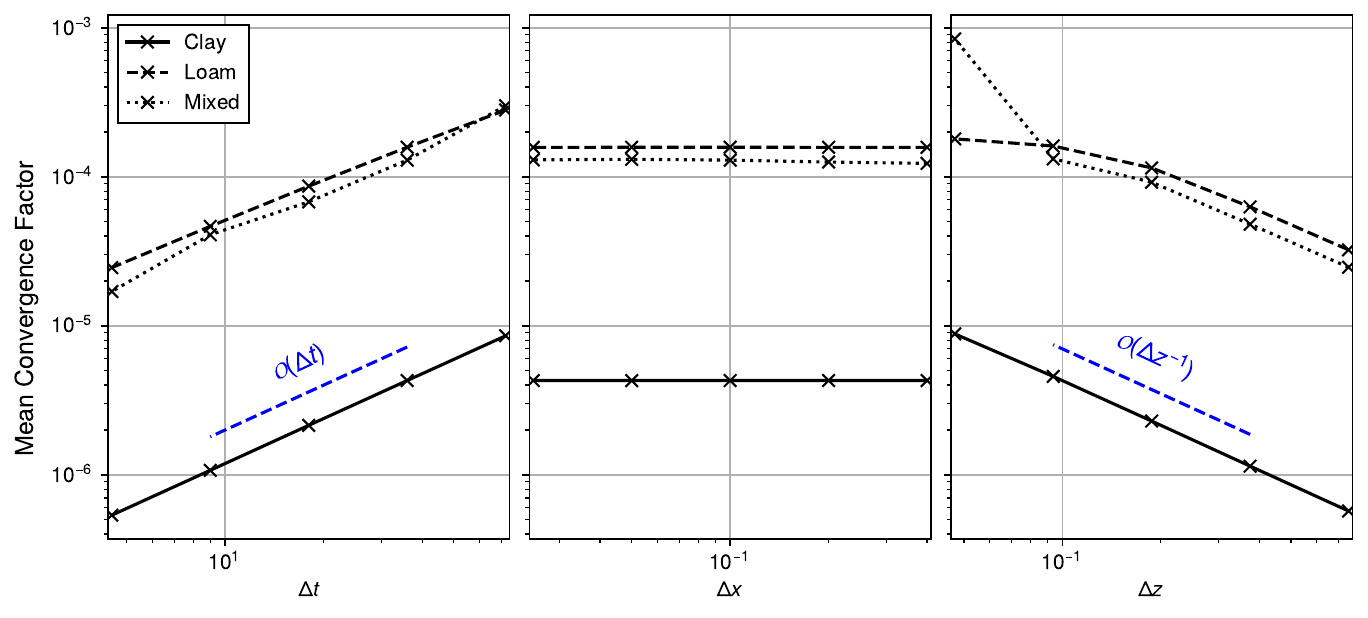}
    \caption{Mean observed convergence factors \(\expcvg\) in the drainage trench example for varying grid sizes. \label{fig:drainage_grid_dependence}}
\end{figure}

Regarding simulation parameters, our analysis says that a smaller time step and vertical grid size should yield a smaller convergence factor, i.e., faster convergence (\Cref{fig:factor_varying_grids}).
Further, an underlying assumption of the linear model was that horizontal processes, and thus also the horizontal resolution \(\dx\), do not affect the convergence speed.
We tested these three hypotheses by varying the grid sizes \(\dt, \dx, \dz\) in isolation, relative to the original simulation parameters.
The results are displayed in \Cref{fig:drainage_grid_dependence}, where we plot the time-averaged observed convergence factors \(\expcvg\).

The leftmost panel shows that the convergence rate indeed decreases with smaller $\Delta t$, irrespective of the soil configuration.
Closer investigation shows that the observed convergence factor decays linearly with \(\dt\).
The analysis results in \Cref{fig:factor_varying_grids}, on the other hand, predict that \(|S|\propto \sqrt{\dt}\).

Changing the horizontal resolution \(\dx\) does not meaningfully affect the convergence speed, indicating that the coupling is indeed driven by vertical dynamics of the Richards equation.
This justifies reducing the dimensionality of the problem to a single vertical column. 

The convergence behavior for the vertical resolution does not match our analysis: 
The convergence factor \(\expcvg\) \emph{increases} with \emph{decreasing} \(\dz\).
This increase is proportional to \(1/\dz\) in all three cases.
For the loam and heterogeneous soil case, the growth slows down with decreasing \(\dz\), except for the leftmost data point in the mixed soil setup.
There, the mean convergence factor jumps to \num{8e-4}.
This is due to a single time step with a larger convergence rate (\(\expcvg_{67} = 0.21\)).
Omitting this value before computing the time average gives \(\expcvg=\num{1.5e-4}\), in line with the other data points.

The linear 1D-0D analysis is not able to explain the observed relationship between \(\dt\), \(\dz\), and \(\expcvg\) in the 2D-1D experiments.
Since the horizontal resolution does not seem to interfere with the convergence behavior, we are fairly certain that this difference stems from nonlinear effects, not from the change in dimensionality.
However, with the tools at hand it is unclear what exactly is causing this discrepancy.

\subsubsection{Large-Scale Experiment: Hillslope}

\begin{figure}
    \centering
    \begin{tikzpicture}[circuit ee IEC]
    \tikzstyle{doublefullarrow}=[<->, >=triangle 60, very thin]
    \tikzset{>=latex}

    \node[] (A) at (0, 0) {};
    \node[] (B) at (0, 1.2) {};
    \node[] (C) at (8, 2) {};
    \node[] (D) at (8, 0) {};
    \node[] (E) at (0, 0.9) {};
    \node[] (F) at (8, 1.7) {};
    
    \draw[draw=black, yslant=0.1] (A) rectangle ++(C);
    
    \node[black] at (6, 1.2) {$\Omega$};
    \node[anchor=south] (dOmega) at (8.5, 1.5) {\footnotesize $\partial\Omega$};
    \draw[black, -] (8, 1) -- (dOmega);

    \draw[waterColor, very thick] ($(B) + (0, 0.03)$) -- ($(C) + (0, 0.03)$);
    \node[waterColor, anchor=south] (Gamma) at (4, 1.6) {\footnotesize $\Gamma$};

    \node[shape=coordinate] (CSorigin) at ($(A) + (0.2, 0.2)$) {};
    \node[] (CSz) at ($(CSorigin) + (0, 0.55)$) {\footnotesize $z$};
    \node[] (CSx) at ($(CSorigin) + (0.55, 0)$) {\footnotesize $x$};
    \draw[black, ->] (CSorigin) -- (CSz);
    \draw[black, ->] (CSorigin) -- (CSx);

    \dimline[extension start length=0.85, extension end length=0.85, label style={above=0.02ex,},] {($(A) + (-1, 0)$)} {($(B) + (-1, 0)$)} {\footnotesize \qty{5}{m}};
    \dimline[extension start length=0.35, extension end length=0.35, label style={above=0.02ex,},] {($(A) + (-0.3, 0)$)} {($(E) + (-0.3, 0)$)} {\footnotesize \qty{4}{m}};
    \dimline[extension start length=-0.08, extension end length=-0.18, label style={above=0.02ex,},] {($(A) + (0, -0.6)$)} {($(D) + (0, -0.6)$)} {\footnotesize \qty{400}{m}};
    \dimline[extension start length=-12, extension end length=-0.8, label style={above=0.02ex,},] {($(D) + (0.5, 0)$)} {($(D) + (0.5, 0.8)$)} {\scriptsize \qty{0.2}{m}};
    \draw [dashed] (E) -- (F);
\end{tikzpicture}
    \caption{Geometry for the hillslope test case. The dashed line marks the initial groundwater height. \label{fig:hillslope_domain}}
\end{figure}
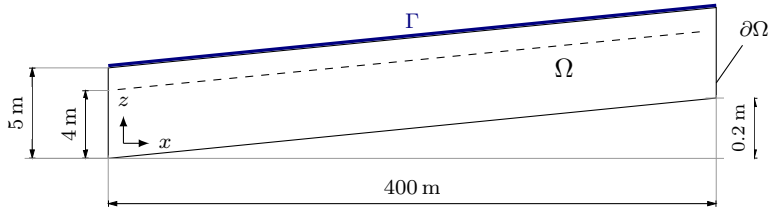

This example is inspired by the hillslope benchmark cases presented in \textcite{maxwell.etal_2014} for the large-scale overland flow model \eqref{eq:maxwell}.
Here we simulate rainfall and resulting outflow onto a sloping plane; its geometry is depicted in \Cref{fig:hillslope_domain}.
The initial and boundary conditions for this experiment are given by:
\begin{subequations}
    \begin{align*}
        \psi(0,\boldsymbol{x}) &= 4 - \left(z-\frac{0.2}{400} x\right) \quad \text{on }\Omega,\\
        \boldsymbol{v} \cdot \boldsymbol{n} &= 0 \quad \text{on } \partial\Omega \setminus \Gamma, \\
        h(0, x) &= 0 \quad \text{on } \Gamma, \\
        \partial_x(hu)(t, 0) &= 0, \\
        (hu)(t, 400) &= 0.
    \end{align*}
\end{subequations}
We pick the initial condition to be linear in \(z\) with the initial groundwater level increasing with \(x\), analogously to \textcite{maxwell.etal_2014}.
Again, we introduce an external source term for the rainfall rate, given by
\begin{equation*}
    r(t) = \begin{cases}
        r & t \leq \qty{200}{min}\\
        0 & t > \qty{200}{min}.
    \end{cases}
\end{equation*}
In all experiments $T=\qty{300}{min}$ and the relaxation parameter $\omega=1$.
The remaining simulation parameters are given in \Cref{tab:hillslope_params}.

\begin{table}[]
    \begin{tabularx}{\textwidth}{| l  c  | l | X | X |}
    \hline
    \textbf{Parameter} & \textbf{} & \textbf{Unit} & \textbf{Sandy Loam} & \textbf{Silt Loam} \\ \hline
    Time step size & $\dt$ & $\unit{s}$ & $60$ & $1$ \\
    Horizontal mesh size & $\dx$ & $\unit{m}$ & $80$ & $80$ \\
    Vertical mesh size & $\dz$ & $\unit{m}$ & $0.2$ & $0.2$ \\
    Manning's coefficient & $n_M$ & $\unit{m^{1/3} min}$ & $\num{3.31e-3}$  & $\num{3.31e-3}$  \\
    Bed slope  & $S_f$  &\% & $0.05$ & $0.05$ \\
    Rainfall rate & $r$ & $\unit{m \per min}$ & $\num{3.30e-4}$ & $\num{3.30e-4}$,\newline $\num{3.30e-5}$ \\
    \hline
    \end{tabularx}
    \caption{Parameters for the coupled hillslope experiments.}
    \label{tab:hillslope_params}
\end{table}

\begin{figure}
    \centering
    \begin{subfigure}{0.47\textwidth}
        \includegraphics[width=\textwidth]{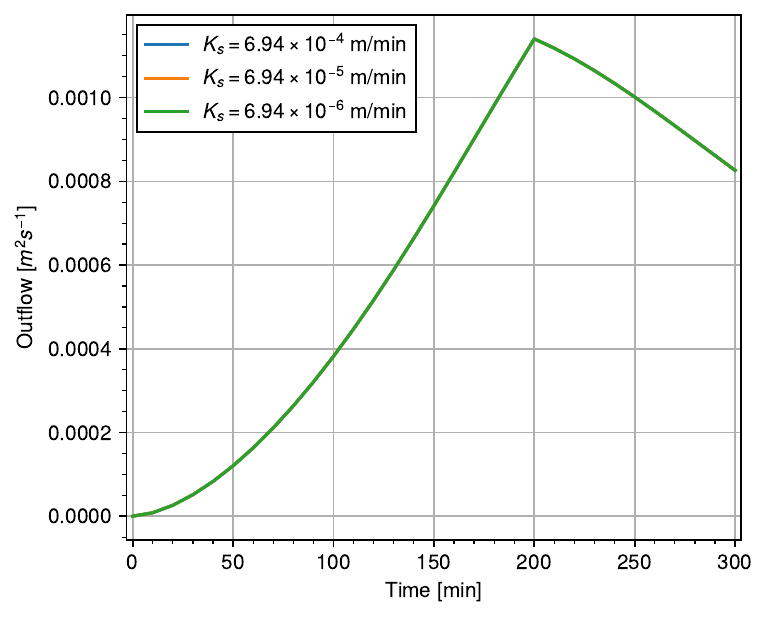}
        \caption{Measured outflow, sandy loam.}
        \label{fig:sandy_loam_outflow}
    \end{subfigure}
    \hfill
    \begin{subfigure}{0.47\textwidth}
        \includegraphics[width=\textwidth]{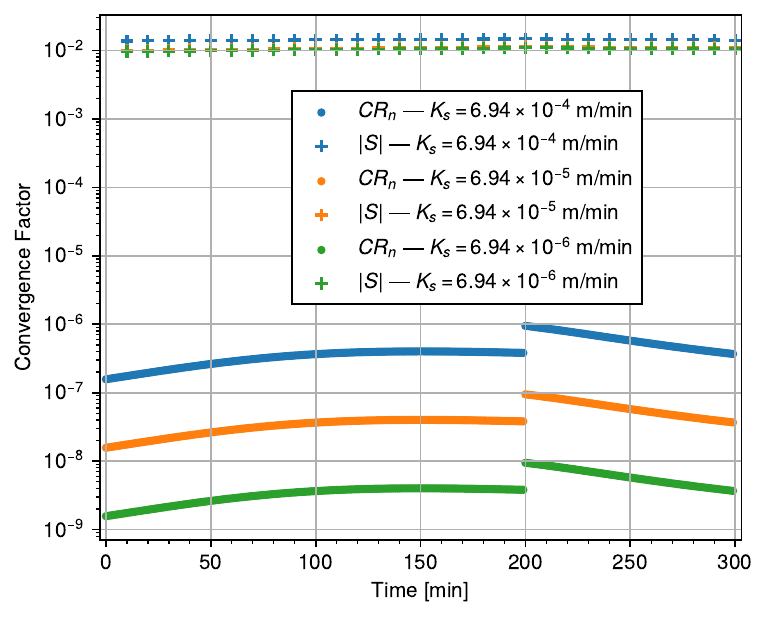}
        \caption{Convergence factors, sandy loam.}
        \label{fig:sandy_loam_cf}
    \end{subfigure}
    \hfill
    \begin{subfigure}{0.47\textwidth}
        \includegraphics[width=\textwidth]{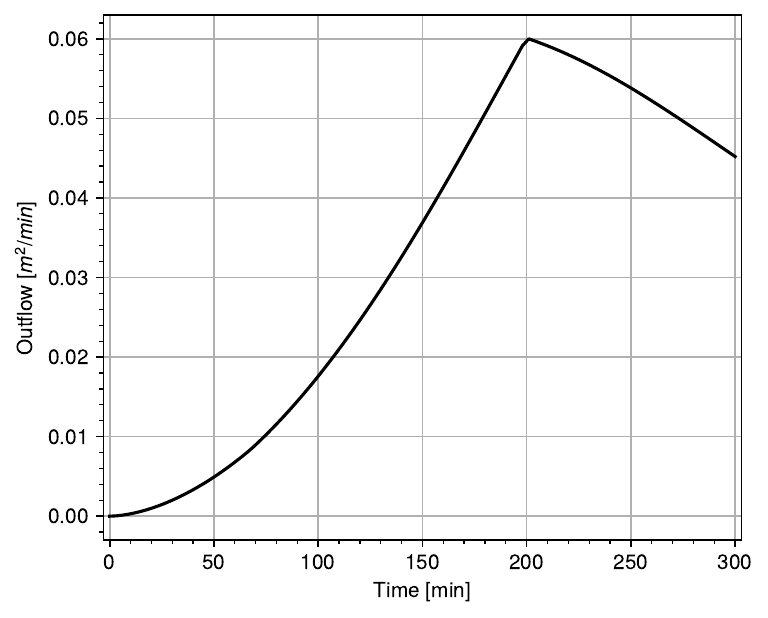}
        \caption{Silt loam, \(r=\qty{3.3e-4}{m/min}\).}
        \label{fig:silt_loam_outflow}
    \end{subfigure}
    \hfill
    \begin{subfigure}{0.47\textwidth}
        \includegraphics[width=\textwidth]{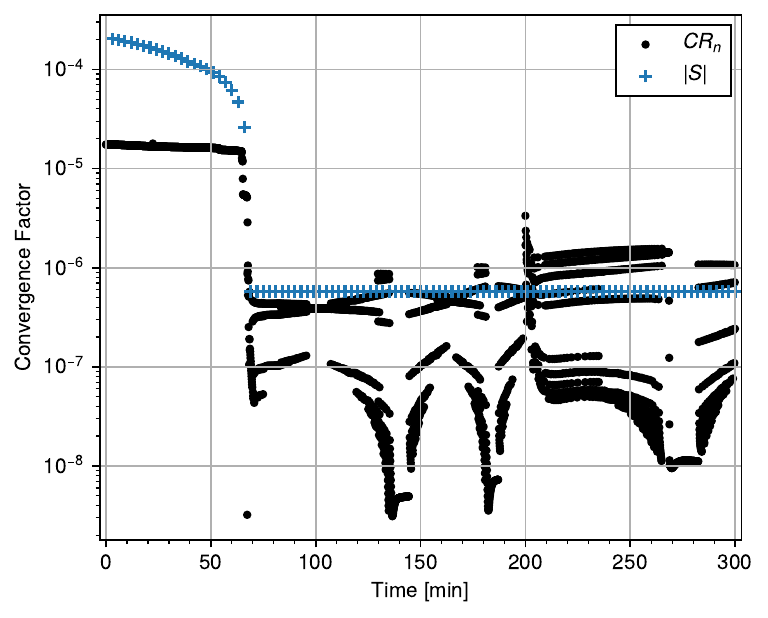}
        \caption{Silt loam, \(r=\qty{3.3e-4}{m/min}\).}
        \label{fig:silt_loam_cf}
    \end{subfigure}
    \hfill
    \begin{subfigure}{0.47\textwidth}
        \includegraphics[width=\textwidth]{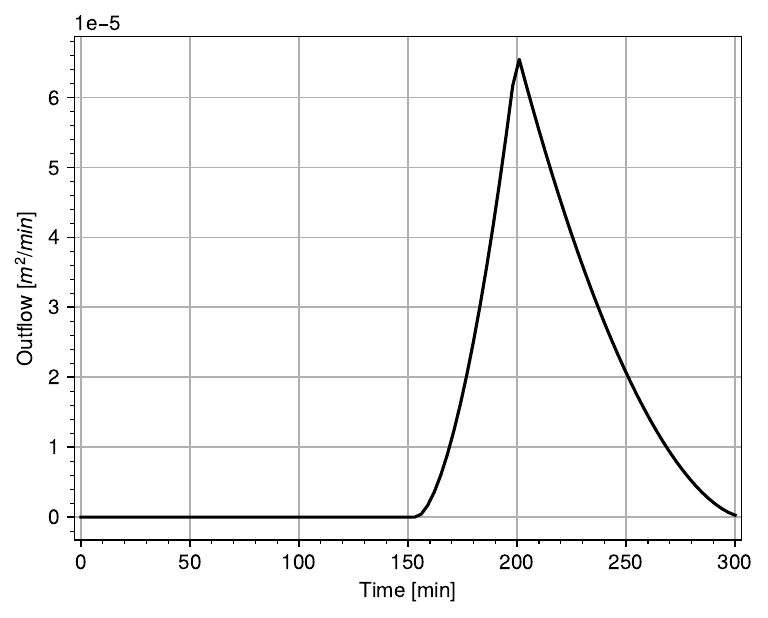}
        \caption{Silt loam, \(r=\qty{3.3e-5}{m/min}\).}
        \label{fig:silt_loam_low_qr_outflow}
    \end{subfigure}
    \hfill
    \begin{subfigure}{0.47\textwidth}
        \includegraphics[width=\textwidth]{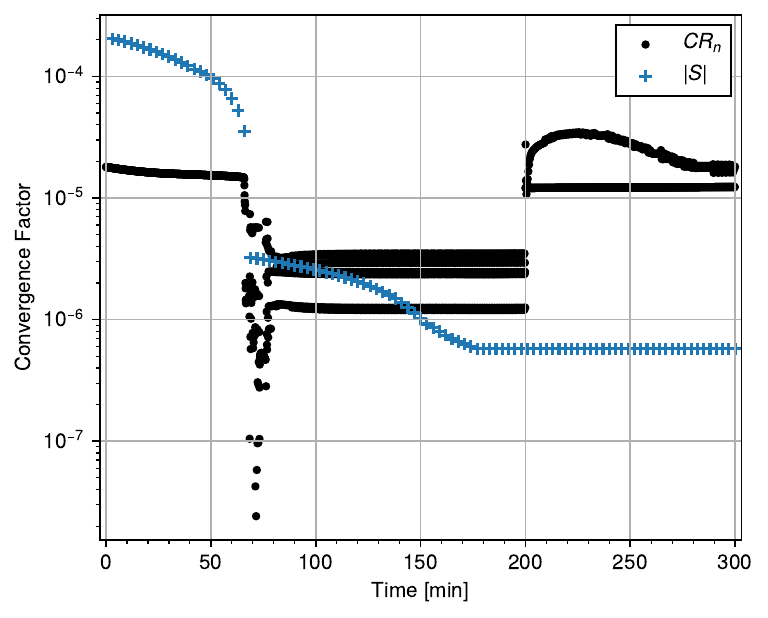}
        \caption{Silt loam, \(r=\qty{3.3e-5}{m/min}\).}
        \label{fig:silt_loam_low_qr_cf}
    \end{subfigure}
    \hfill
    \caption{Numerical results for the coupled hillslope cases.}
\end{figure}

For our first experiment, we use the sandy loam material parameters from \Cref{tab:van_genuchten_params} to consider three different values of $K_s$ and its impact on the convergence factor.
\Cref{fig:sandy_loam_outflow} shows the measured outflow at the left boundary over time, \(q_\mathrm{out}(t) = h(t,0)u(t,0)\).

In our model, the saturated hydraulic conductivity \(K_s\) does not produce visible differences in the measured outflow, as opposed to \textcite[Fig. 4]{maxwell.etal_2014}.
We think that this is due to differences between our model and theirs.
As explained in \Cref{sec:governing_equations} we do not include boundary condition switching in the interface boundary condition.
This way, we cannot expect our subsurface model implementation to show the same infiltration behavior as more sophisticated coupled surface-subsurface flow models.

Nevertheless, we do see differences in the convergence factors due to changes in \(K_s\), as shown in \Cref{fig:sandy_loam_cf}:
Both \(|S|\) and \(\expcvg_n\) increase with increasing $K_s$.
The observed convergence factor \(\expcvg_n\) is proportional to \(K_s\).
On the other hand, \(|S|\) increases only slightly even though \(K_s\) increases by a factor of \num{10}.
This could again be due to the fact that we use the spatial means of \(c\) and \(K\) for the computation of \(|S|\), while \(\expcvg_n\) might mainly be driven by near-interface values where \(\psi>0\) and thus \(K(\psi)=K_s\), cf. \Cref{eq:genuchten}.

In general, we see larger differences here between \(\expcvg_n\) and the predicted convergence factor \(|S|\), compared to the drainage trench experiments before:
\Cref{eq:discr_analysis_result} predicts convergence factors around \num{e-2}, while we measure convergence rates between \num{e-9} and \num{e-6}.

Additionally, we note that \(\expcvg_n\) jumps at the discontinuity in the rainfall rate \(r(\qty{200}{min})\), a behavior not captured by the predicted convergence factor.
In the linear analysis, adding rainfall does not affect the computation of the convergence rate: \(r(t)\) does not change from one iteration to the next, thus being irrelevant for computing \(S\) or \(\Sigma(\omega)\).
This is not the same in the nonlinear case.

To investigate the effect of the rainfall rate on the convergence behavior in more detail, we carried out another experiment for two different values of \(r\).
Due to stiffness issues of our subsurface solver, we carried out this experiment  with the silt loam parameters used in the drainage trench case instead of sandy loam, cf. \Cref{tab:van_genuchten_params}.
See \Cref{tab:hillslope_params} for the simulation parameters in this case.

\Cref{fig:silt_loam_outflow,fig:silt_loam_low_qr_outflow} show the measured outflow \(q_\mathrm{out}(t)\) at \(x=0\) over time for \(r=\qty{3.3e-4}{m/min}\) and \(r=\qty{3.3e-5}{m/min}\), respectively.
The lower rainfall rate leads to a much smaller outflow rate and different shape of the curve, since the ground is able to take up and release water more quickly compared to the influx.
The peak in outflow is reached at \(t=\qty{200}{min}\) when the rainfall is turned off.

The expected and observed convergence behavior is displayed in \Cref{fig:silt_loam_cf,fig:silt_loam_low_qr_cf}.
The results show highly nonlinear dynamics of \(\expcvg_n\).
It spans multiple orders of magnitude during the simulation and shows large jumps from one time step to the next in both experiments.
Again, a jump in the convergence factor is visible at \(t=\qty{200}{min}\).

Compared to the sandy loam experiments above, our estimator \(|S|\) does a better job capturing the order of magnitude and dynamics of \(\expcvg_n\).
In particular, there is a strong decay in \(\expcvg_n\) at around \(t=\qty{60}{min}\) in both experiments that is also captured by \(|S|\).
Since \(\dt, \dz\), and \(L_z\) are constant, the decay in \(|S|\) must be due to changes in \(c\) and \(K\) due to increasing soil saturation.
It appears that significant changes in material properties drive the convergence rate of the coupling iteration here, an effect our analysis is able to capture.

In the case with the larger rainfall rate, \Cref{fig:silt_loam_cf}, \(|S|\) stays roughly constant at around \num{5e-7} for \(t>\qty{60}{min}\), whereas \(\expcvg_n\) continues to vary between \num{e-9} and \num{5e-6}.
For \(r=\qty{3.3e-5}{m/min}\) in \Cref{fig:silt_loam_low_qr_cf}, the qualitative behavior is more different:
The observed \(\expcvg_n\) is in two different value ranges in the time intervals \([\qty{60}{min}, \qty{200}{min}]\) and \([\qty{200}{min}, \qty{300}{min}]\).
In both intervals, the values jump from one time step to the next but do not show a time-dependent trend.
The predicted convergence factor \(|S|\) on the other hand decays smoothly for \(t\in[\qty{60}{min}, \qty{180}{min}]\), before stabilizing around \num{5e-7}.
Turning off the rainfall rate at \(t=\qty{200}{min}\) has no visible effect on \(|S|\) in both experiments.

As for the drainage trench example, the observed material parameters in the coupled hillslope experiment are small, even when disregarding the spatial mean: 
The maximal local values for \(c,K\) across all simulations are \(c=10.7\), \(K=\num{1.16e-5}\) (both observed for sandy loam), corresponding to a small expected convergence rate (cf. \Cref{fig:factor_analysis}).

\section{Summary and Conclusions}
\label{sec:conclusion}

This paper considers coupled surface-subsurface models which solve the Richards and shallow water equations coupled via boundary conditions.
We employ a partitioned setup where the two subsolvers are coupled with a sequential procedure, allowing for acceleration of the iterations via relaxation.
To our knowledge, this is the first time that the discrete and continuous convergence analysis techniques for coupled PDE systems have been applied to this problem.
We are contributing to the field by examining what conclusions we can draw in practice when applying existing linear analysis techniques to this highly nonlinear problem.

We derived a linearized, 1D-0D model problem, assuming that vertical processes are most relevant.
Linear analysis techniques from \textcite{gander.etal_2016,monge.birken_2018} allowed us to determine convergence factors and an optimal relaxation parameter for the fully continuous and fully discrete versions, and to show their dependence on the material properties of the subsurface, but also on the resolution in time and space. In particular, the convergence rate will be small when \(c\) or \(K\) are small. These theoretical results are confirmed numerically for the 1D-0D problem. 

Nonlinear 2D-1D experiments for different realistic materials show that the convergence speed is indeed mainly driven by vertical processes and resolution in Richards' equation.
Quantitatively, the linear analysis tends to overestimate the observed convergence rate and naturally fails to represent processes which are not part of the linear model (e.g., the effect of external source terms and spatial variability).
Also, in the drainage trench example, the convergence factor decreased linearly with respect to the time step size while the analysis predicted a square-root dependence.
Even more surprisingly, it increased when increasing the vertical resolution.

Overall the observed convergence factors do not exceed \num{e-3} for the materials considered and typical simulation parameters, making acceleration using relaxation unnecessary. 
The linear analysis is able to explain the small observed convergence factors, as well as how dynamically changing material properties during a simulation affect the speed of convergence. 

\section*{Acknowledgements}
The research presented in this paper is a contribution to the Strategic Research Area “ModElling the Regional and Global Earth system”, MERGE, funded by the Swedish government.
We wish to thank Robert Klöfkorn and Niklas Kotarsky for their contributions to the code used in this paper, as well as Martin Gander for input on early results, especially regarding the continuous analysis.

\printbibliography{}

\appendix

\section{Well-Posedness of the Coupled Problem}\label{subsec:well_posedness}

We study the linearized 1D-0D problem without iteration or relaxation:
\begin{subequations}
\begin{equation}
    \begin{aligned}
        c \frac{\partial\psi}{\partial t} - K\frac{\partial^2}{\partial z^2} \psi &= 0 \quad\text{on }(0,T]\times(-L,0),\\
        \psi(t,-L) &= 0 \quad\text{on }[0,T],\\
        \psi(t,0) &= h(t) \quad\text{on }[0,T],\\
        \psi(0,z) &= \psi_0(z) \quad\text{on }[-L,0],
    \end{aligned}
\end{equation}
\begin{equation}
    \begin{aligned}                        
        \frac{dh}{dt} & = -K \left(\frac{\partial}{\partial z} \psi(t,0) + 1\right) \quad\text{on }(0,T],\\
        h(0) &= h_0.
    \end{aligned}
\end{equation}
\end{subequations}
Applying the Laplace transform \eqref{eq:laplace_transform}, we arrive at
\begin{subequations}\label{eq:coupled_laplace}
\begin{equation}
    \begin{aligned}\label{eq:richards_laplace}
        cs\hat{\psi} - K \partial_z^2 \hat{\psi} &= 0, \\
        \hat{\psi}(s,-L) &= 0,
    \end{aligned}
\end{equation}
\begin{equation}
    s\hat{h} = -K \left(\frac{\partial}{\partial z} \hat{\psi}(s,0) + 1\right) \label{eq:swe_laplace}.
\end{equation}
\end{subequations}

If we find a solution to this system to which the backward transform is well-defined, the problem is well-posed.
The roots to the characteristic equation for \Cref{eq:richards_laplace} are \(\lambda_{1/2} = \pm \sqrt{cs/K}\), yielding a general solution:
\begin{equation*}
    \hat{\psi} 
    = A \sinh\left(\sqrt{\frac{cs}{K}}z\right) 
    + B \cosh\left(\sqrt{\frac{cs}{K}}z\right).
\end{equation*}
The Dirichlet boundary condition \(\hat{\psi}(s,-L) = 0\) gives
\begin{equation*}
    \hat{\psi}(s,z) = \frac{B}{\sinh\left(\sqrt{\frac{cs}{K}}L \right)} \sinh\left(\sqrt{\frac{cs}{K}} (z+L)\right),
\end{equation*}
and from the interface boundary condition \(\hat{\psi}(s,0)=\hat{h}\) we obtain \(B=\hat{h}\).
Inserting this result into \Cref{eq:swe_laplace} gives
\begin{equation*}
    s\hat{h} = - \hat{h} \sqrt{cKs} \coth\left(\sqrt{\frac{cs}{K}}L\right) - K,
\end{equation*}
and thus
\begin{equation*}
    \hat{h} = \frac{-K}{s+\sqrt{\frac{cK}{s}} \coth\left(\sqrt{\frac{cs}{K}}L\right)}.
\end{equation*}

Assuming that \(c,K,L>0\), the denominator is well-defined for \(s\neq 0\) and thus the solution to the Laplace-transformed coupled problem has its only pole at \(s=0\).
The inverse Laplace transform is therefore well-defined for \(s\in\mathbb{C}^+\), giving existence and uniqueness for \(\psi(t,z)\) and \(h(t)\).

\section{Convergence factor estimate in code}\label{app:convergence_estimate}

We estimate convergence according to \Cref{eq:conv_estimate} for the linearized 1D-0D coupling problem.
Here, \(\|h^{k,n}\|_2 \equiv |h^{k,n}|\) and thus
\begin{equation*}
    \expcvg_n = \frac{|\tilde{h}^{k,n}-h^{k-1,n}|}{|\tilde{h}^{k-1,n} - h^{k-2,n}|}.
\end{equation*}
Recall from \Cref{eq:update_height_norelax,eq:update_height_relax}:
\begin{equation}\label{eq:update_steps}
    \begin{gathered}
        \tilde{h}^{k,n} = Sh^{k-1,n} + \xi^{n-1},\\
        h^{k,n} = \omega \tilde{h}^{k,n} + (1-\omega)h^{k-1,n} = \underbrace{(1 - \omega + \omega S)}_{\eqqcolon\Sigma(\omega)} h^{k-1,n} + \omega \xi^{n-1}.
    \end{gathered}
\end{equation}

We now define the errors at time \(t^n\) with respect to the true solution \(h(t)\):
\begin{equation}
    e^{k,n} \coloneqq h^{k,n} - h(t^n), \quad \tilde{e}^{k,n} \coloneqq \tilde{h}^{k,n} - h(t^n)
\end{equation}
The above iteration \eqref{eq:update_steps} can be rewritten in terms of \(e^{k,n}\) and \(\tilde{e}^{k,n}\):
\begin{equation*}
    e^{k,n} + h(t^n) = \omega (\tilde{e}^{k,n} + h(t^n)) + (1-\omega) (e^{k-1,n} + h(t^n)).
\end{equation*}
That is,
\begin{equation*}
    e^{k,n} = \omega \tilde{e}^{k,n} + (1-\omega) e^{k-1,n}.
\end{equation*}
Similarly,
\begin{equation*}
    \tilde{e}^{k,n} + h(t^n) = S (\tilde{e}^{k-1,n} + h(t^n)) + \xi^{n-1}.
\end{equation*}

Assume now that \(h^{k-1,n} = h(t^n)\).
If we have defined our coupling iteration appropriately, it must hold that \(\tilde{h}^{k,n} = h(t^n)\) as well.
Thus,
\begin{equation*}
    h(t^n) = S h(t^n) + \xi^{n-1}
\end{equation*}
and therefore
\begin{equation*}
    \tilde{e}^{k,n} = S e^{k-1,n} \quad\Rightarrow \quad e^{k,n} = \Sigma(\omega) e^{k-1,n}.
\end{equation*}
We can then rewrite
\begin{equation*}
    |\tilde{h}^{k,n}-h^{k-1,n}| = |(\tilde{e}^{k,n} + h(t^n)) - (e^{k-1,n} + h(t^n))| = |(S-1)e^{k-1,n}|
\end{equation*}
and obtain the desired result \eqref{eq:cr_motivation}
\begin{equation*}
    \expcvg_n = \frac{|\tilde{h}^{k,n}-h^{k-1,n}|}{|\tilde{h}^{k-1,n} - h^{k-2,n}|} = \frac{|(S-1)e^{k-1,n}|}{|(S-1)e^{k-2,n}|} = \frac{|e^{k-1,n}|}{|e^{k-2,n}|} = |\Sigma(\omega)|
\end{equation*}
up to floating point precision.

Note: The same result holds if one defines the residual in terms of \emph{both} the surface flux \(v_s^{k,n}\) and the water height, i.e.,
\begin{equation*}
    \mathrm{res}^{n,k} = \left\| \begin{pmatrix}\hat{h}^{n,k}\\\hat{v}_s^{n,k}\end{pmatrix} - \begin{pmatrix}h^{n,k-1}\\v_s^{n,k-1}\end{pmatrix} \right\|. 
\end{equation*}

\end{document}